\documentclass[12pt,oneside]{amsart}
\usepackage{latexsym}
\usepackage{amsmath,amssymb,amsthm}
\usepackage{amscd}
\usepackage[dvips]{graphics}
\newcommand{\rp}{\mathbb{RP}}
\newcommand{\re}{\mathbb{R}}
\newcommand{\co}{\mathbb{C}}

\newcommand{\pgl}[1]{\mathbf{PGL}(#1,\mathbb{R})}
\newcommand{\Sl}[1]{\mathbf{SL}(#1,\mathbb{R})}
\newcommand{\gl}[1]{\mathbf{GL}(#1,\mathbb{R})}
\newcommand{\na}{\nabla}
\newcommand{\sfrac}[2]{{\textstyle \frac{#1}{#2}}}
\newcommand{\D}{\displaystyle}

\newtheorem{prop}{Proposition} 
\newtheorem{cor}[prop]{Corollary}
\newtheorem{thm}{Theorem}

\newtheorem{lem}[prop]{Lemma}

\theoremstyle{remark}
\newtheorem*{rem}{Remark}
\newtheorem*{ack}{Acknowledgements}

\begin{document}
\author{John C.\ Loftin}
\title{The Compactification of the Moduli Space of Convex $\rp^2$ Surfaces, I}
\maketitle


\section{Introduction}

In \cite{goldman90a} Goldman proves that the deformation space
$\mathcal G_g$ of convex $\rp^2$ structures on an oriented closed
surface of genus $g\ge2$ is a cell of $16g-16$ real dimensions.
He constructs explicit coordinates of the space based on a
Fenchel-Nielsen type pants decomposition of the surface.  In
particular, for each boundary geodesic loop in this decomposition,
there is a holonomy map in $\pgl{3}$, unique up to conjugacy,
which measures how a choice of $\rp^2$ coordinates develops around
the loop.

In \cite{labourie97,loftin01} there is another proof of Goldman's
theorem which introduces new coordinates on $\mathcal G$.  Using a
developing map due to C.P.\ Wang \cite{wang91} and deep results in
affine differential geometry of Cheng-Yau
\cite{cheng-yau77,cheng-yau86}, there is a natural correspondence
between a convex $\rp^2$ structure on a surface of genus $g\ge2$
and a pair $(\Sigma,U)$ of a conformal structure $\Sigma$ and a
holomorphic cubic differential $U$ on $\Sigma$.  This shows that
$\mathcal G_g$ is the total space of a $5g-5$ complex dimensional
vector bundle over Teichm\"uller space, and so by Riemann-Roch,
$\mathcal G_g$ is a cell of complex dimension $8g-8$.

Wang's formulation does give a way to calculate the holonomy by a
first-order linear system of PDEs, but unfortunately the data
depends on the solution to a separate nonlinear PDE
(\ref{wang-eq}) below, and thus it seems we cannot hope to use
this to relate the two coordinate systems.  However, on a
noncompact Riemann surface of finite type, we can force solutions
of (\ref{wang-eq}) to behave well near the punctures, and in fact
calculate the holonomy around a loop at each puncture.

\begin{thm}
Given a Riemann surface $\Sigma$ of finite type, and a holomorphic
cubic differential $U$ on $\Sigma$ with poles of order 3 allowed
at each puncture. Let $R_i$ be the residue at each puncture. Then
there is an $\rp^2$ structure corresponding to $(\Sigma,U)$.  The
$\rp^2$ holonomy type around each puncture is determined by the
$R_i$.
\end{thm}

A more explicit version of this theorem is given below as Theorem
\ref{holonomy-thm}.  The space of all residues $R\in\co$ maps
two-to-one onto the space of holonomy types if ${\rm Re}\,R\neq0$.
These two cases are distinguished by another of Goldman's
coordinates, the vertical twist parameter, which in part
corresponds to the twist parameter in Fenchel-Nielsen theory.  The
vertical twist parameter approaching $\infty$ roughly corresponds
stretching the $\rp^2$ structure along a neck.

\begin{thm}
Given a pair $(\Sigma,U)$ as above, at a puncture with residue
$R_i$ with ${\rm Re}\,R_i\neq0$, the vertical twist parameter of
the $\rp^2$ structure at that puncture is $\pm\infty$, the sign
being equal to the sign of ${\rm Re}\,R_i$.
\end{thm}

This theorem combines Propositions \ref{plus-prop} and
\ref{minus-prop} below.

There are natural examples of such pairs $(\Sigma,U)$ at the
boundary of the total space of the space of cubic differentials
over the moduli space of Riemann surfaces $\overline{\mathcal
M}_g$.  In particular, at a nodal curve at the boundary of
$\mathcal M_g$, regular cubic differentials have poles of order at
most 3 at the nodes, with opposite residues across each side of
the node.

\begin{thm} On the total space of the (V-manifold) vector bundle of
regular cubic differentials over $\overline{\mathcal M}_g$, the
holonomy type and vertical twist parameters vary continuously.
\end{thm}

This theorem is made more precise below in Theorem \ref{families}.
In particular, we must also make a mild technical assumption on
families of cubic differentials which approach regular cubic
differentials with residue 0.

The proof of these results relies on a fundamental correspondence
due to Cheng-Yau \cite{cheng-yau77,cheng-yau86}, that each convex
bounded domain $\Omega\subset \rp^n$ may be canonically identified
with a hypersurface, the hyperbolic affine sphere $H\subset
\re^{n+1}$ which is asymptotic to the boundary of the cone
$\mathcal{C}(\Omega)$ over $\Omega$.  Any projective group
$\Gamma\subset\pgl{n+1}$ which acts on $\Omega$ lifts to a linear
group which acts on $\re^{n+1}$ and takes $H$ to itself, and
natural differential geometric structure descends to any manifold
$\Omega/\Gamma\equiv H/\Gamma$.  This theory is reviewed in
Section \ref{h-a-s} below.  See also \cite{labourie97,loftin01}.

In dimension $n=2$, there is a natural system of first-order PDEs
which develop any hyperbolic affine sphere in $\re^3$
\cite[Wang]{wang91}. (Wang's theory
is analogous to the more familiar case of minimal
surfaces in $\re^3$.)  We give a version of this formulation in
Section \ref{wang-corr} below.  In particular, given a Riemann
surface $\Sigma$ with conformal metric $h$, a cubic differential
$U$ over $\Sigma$, and a function $u$ so that
\begin{equation} \label{wang-eq-intro}
\Delta u=-4 e^{-2u} \|U\|^2 + 2 e^u + 2\kappa,
\end{equation}
there is a system of first-order PDEs
(\ref{z-deriv}-\ref{zbar-deriv}) below which integrate the
structure equations for a hyperbolic affine sphere to give a map
from the universal cover of $\Sigma$ to a hyperbolic affine sphere
$H$. Moreover, Cheng-Yau \cite{cheng-yau86} guarantees that if the
metric $e^uh$ is complete, then $H$ is asymptotic to the boundary
of a convex cone $\mathcal C(\Omega)$ for a domain $\Omega$ as in
the previous paragraph.  So given $(\Sigma,U)$ and a function $u$
satisfying (\ref{wang-eq-intro}), there is a convex $\rp^2$
structure on $\Sigma$.  Moreover, the first-order PDE system
(\ref{z-deriv}-\ref{zbar-deriv}) naturally calculates the holonomy
of the $\rp^2$ structure around any loop.  This relates the affine
differential geometry of Wang and Cheng-Yau to Goldman's
coordinates on the deformation space $\mathcal G$ of convex
$\rp^2$ structures. Goldman's coordinates on the deformation space
of convex $\rp^2$ structures are discussed in Section
\ref{goldman-coor} below.

Wang solved his PDE (\ref{wang-eq-intro}) to find a unique
solution on any compact Riemann surface of genus $g\ge2$.  In
Section \ref{finding-sol} below, we extend this to a complete
hyperbolic Riemann surface $\Sigma$ of finite type with cubic
differential $U$ which has poles of order 3 at each puncture.  We
find an Ansatz metric $h$ on $\Sigma$ in order to find a solution
$u$ to (\ref{wang-eq-intro}) which behaves well near each
puncture.  In particular, we construct a barrier to bound $u$ near
each puncture.  Then calculating the holonomy around each puncture
involves solving a linear systems of ODEs which is asymptotic to a
constant coefficient system as we approach the puncture.  This
allows us to determine the eigenvalues of the holonomy matrix in
terms of the residue $R$ of $U$ ($R$ is the leading coefficient of
the Laurent series of $U$ at the puncture). The actual conjugacy
type is more subtle if the eigenvalues are repeated.  In
particular, a geometric result of Choi \cite{choi94b}
rules out certain holonomy types. See Section \ref{holo-sec}
below.

In Section \ref{detail-sec} we use Wang's developing map
(\ref{z-deriv}-\ref{zbar-deriv}) to provide more information about
the structure of the $\rp^2$ structure near the punctures. In
particular (if the residue $R\neq0$), along certain paths which
approach the puncture, the (\ref{z-deriv}-\ref{zbar-deriv})
restrict to a system of ODEs of the form
$$ \partial_y \mathbf{X}=[\mathbf{B}+\mathbf{E}(x,y)]\mathbf{X},$$
where $\mathbf{B}$ is a constant matrix and the error term
$\mathbf{E}(x,y)$ decays exponentially as $y\to\infty$.  These
equations have been extensively studied, originally by Dunkel
\cite{dunkel}.  In Appendix \ref{app}, we extend estimates on solutions of
such systems to the case of a parameter $x$.  In particular, we
can determine that the vertical twist parameter at each such end
is $\pm\infty$, the sign agreeing with that of ${\rm Re}\,R$.

Finally, in Section \ref{deg-fam-sec}, we extend these results to
families of pairs $(\Sigma_t,U_t)$ in which $\Sigma_t$ tends to a
nodal curve at the boundary of the Deligne-Mumford
compactification $\overline{\mathcal M}_g$ of the moduli space of
Riemann surfaces, and the cubic differential $U_t$ over $\Sigma_t$
tends to a regular cubic differential over the nodal curve.  A
regular cubic differential allows poles of order 3 at either side
of the node, and the residues $R,R'$ at either side of the node
satisfy $R=-R'$.  We first review the analytic theory of
$\overline{\mathcal M}_g$ in terms of plumbing coordinates which
replace each node by a thin neck.  Our barriers extend across
the neck to control a solution $u_t$ to (\ref{wang-eq-intro}).
Thereby, the holonomy and vertical twist parameters along
this end behave well in families.  (Controlling the vertical twist parameters
in families involves some fairly subtle analysis to fix a natural
gauge. Here the extension of Dunkel's
theorem to systems with parameters in Appendix A is essential.)

As mentioned above, Goldman worked out the analog of
Fenchel-Nielsen coordinates on $\mathcal G_g$ \cite{goldman90a}.
In general, one would like to extend the rich theory of
Teichm\"uller space to $\mathcal G_g$.  The symplectic theory is
worked out rather well: There is a natural symplectic form on
$\mathcal G_g$ due to Goldman \cite{goldman84,goldman90b} which
extends the Weil-Petersson symplectic form on Teichm\"uller space.
Kim \cite{kim99} has used Goldman's coordinates to show that
$\mathcal G_g$ is naturally symplectomorphic to $\re^{16g-16}$,
which extends a result of Wolpert \cite{wolpert85} about
Teichm\"uller space.

Another example of this is provided by Hitchin \cite{hitchin92},
who finds a connected component, the Teichm\"uller component, of
the space of representations of $\pi_1(\Sigma_g)$ into $\pgl{n}$
for all $n$.  In the case $n=3$, the Teichm\"uller component is
given by the space $H^0(\Sigma,K^2)\oplus H^0(\Sigma,K^3)$ for a
fixed complex structure on $\Sigma$.  This space may be identified
by $\mathcal G_g$ by the holonomy representation of the $\rp^2$
structure \cite[Choi-Goldman]{choi-goldman93}.  In particular,
Hitchin extends the theory of Higgs bundles he earlier used to
study Teichm\"uller space itself in \cite{hitchin87} to study
representation spaces into Lie groups.  We should remark that
Fran\c{c}ois Labourie has shown in unpublished work that if the
$H^0(\Sigma,K^2)$ part of Hitchin's correspondence is zero, then
Hitchin's cubic differential may be identified Wang's cubic
differential.

We should also mention that Darvishzadeh-Goldman \cite{d-goldman}
have constructed an analog of the Weil-Petersson metric on
$\mathcal G_g$ which interacts with the symplectic form to form an
almost complex, and indeed, and almost K\"ahler structure.  It is
not clear if this almost complex structure is integrable or how it
might relate to the complex structure provided by the pairs
$(\Sigma,U)$.

On the compactified moduli space of Riemann surfaces
$\overline{\mathcal M}_g$, there has been extensive work relating
the complex structure of $\overline{\mathcal M}_g$ as a V-manifold
with the geometry provided by the hyperbolic metric (i.e.\ the
Fenchel-Nielsen coordinates).  See in particular Wolf
\cite{wolf91} and  Wolpert \cite{wolpert90}.  The present work is
an attempt to extend this theory to the moduli space of $\rp^2$
structures.  In particular, the space of pairs $(\Sigma,U)$ of a
Riemann surface $\Sigma$ of genus $g\ge2$ and a cubic differential
$U$ over $\Sigma$ is naturally partially compactified to be the
total space of a (V-manifold) vector bundle $\mathcal S \to
\overline{\mathcal M}_g$, where the fiber of $\mathcal S$ is the
space of regular cubic differentials over the possibly singular
curve $\Sigma$. This provides a natural complex structure on a
partial compactification of the moduli space of convex $\rp^2$
structures on oriented surfaces of genus $g$. The $\rp^2$
structure on $\Sigma$ is the analog of the hyperbolic metric, and
we relate this complex structure on $\mathcal S$ to the
coordinates on $\mathcal G_g$ given by Goldman.  In this work we
only address the holonomy and twist parameters around a neck which
is being pinched to a node.  We hope to address the other
coordinates in Goldman's chart in future work.  Also, we only
address the topological relationship between $\mathcal S$ and
Goldman's coordinates in this paper.  Eventually one might be able
to study relationship between the two coordinates
real-analytically.  This point is already complicated in the case
of $\overline{\mathcal M}_g$ \cite[Wolf-Wolpert]{wolf-wolpert92}.

There are nonconvex $\rp^2$ structures on surfaces, which are classified
by Choi \cite{choi94a,choi94b,choi96}.

Also, there is recent work of Benoist on compactifying the space
of faithful representations of a finite-type group
$\rho:\Gamma\to\pgl{n+1}$ so that there is a properly convex
$\Gamma$-invariant domain $\Omega\subset\rp^n$ so that
$\Omega/\rho(\Gamma)$ is compact \cite{benoist03}. This extends
the theorem of Choi-Goldman \cite{choi-goldman93} to higher
dimensions.

\begin{ack}
I would like to thank S.-T.\ Yau, Curt McMullen, Bill Goldman,
Suhyoung Choi, and Scott Wolpert, Michael Thaddeus, Ravi Vakil,
Richard Wentworth for useful conversations.  I would like to thank
Curt McMullen  especially for suggesting that there should be a
relationship between the residues of the cubic differentials and
the $\rp^2$ holonomy of the end.
\end{ack}

\section{Goldman's coordinates} \label{goldman-coor}

An $\rp^n$ structure on a manifold $M$ consists of a maximal atlas
of charts in $\rp^n$ with transition maps in $\pgl{n+1}$.  We may
also call $M$ an $\rp^n$ manifold.  Note that the straight lines
in any coordinate chart in $\rp^n$ are preserved by the transition
maps.  A path in $M$ which in any such coordinate chart is a
straight line is called an $\rp^n$ \emph{geodesic}. Also, given an
$\rp^n$ manifold $M$ and a coordinate chart in $\rp^n$ around a
point $p$, analytic continuation induces a map, the
\emph{developing map} dev extending the coordinate chart map from
the universal cover $\tilde M$ of $M$ (with basepoint $p$) to
$\rp^n$.  Also, there is a \emph{holonomy homomorphism} hol from
$\pi_1M$ to $\pgl{n+1}$ so that $\forall \gamma\in\pi_1M$,
$$ {\rm dev}\circ\gamma = {\rm hol}(\gamma) \circ {\rm dev}.$$
For any other choice of basepoint and/or coordinate chart, there
is a map $g\in\pgl{n+1}$ so that
$$ {\rm dev}' = g\circ {\rm dev} \quad \mbox{and} \quad
{\rm hol}'(\gamma) = g \circ {\rm hol}(\gamma) \circ g^{-1}.
$$ For more details, see e.g.\ Goldman \cite{goldman90a}.  In
particular, the holonomy map is unique up to conjugation in
$\pgl{n+1}$.

$M$ is a \emph{convex} $\rp^n$ manifold if dev$:\tilde M \mapsto
\rp^n$ is a diffeomorphism onto a domain $\Omega$ so that $\Omega$
is a convex subset of some $\re^n\subset \rp^n$. $M$ is called
\emph{properly convex} if in addition $\Omega$ is properly
contained in some $\re^n\subset\rp^n$. In this case, there is a
representation hol of $\Gamma=\pi_1M$ into $\pgl{n+1}$ so that
$\Gamma$ acts discretely and properly discontinuously on $\Omega$.
The quotient $\Omega/\Gamma$ is our $\rp^n$ manifold $M$.  See
e.g.\ \cite{goldman90a} for details.

In \cite{goldman90a}, Goldman introduced coordinates on the
deformation space of convex real projective structures on a given
closed oriented surface $S$ of genus $g\ge2$. These coordinates
are analogous to Fenchel-Nielsen coordinates on Teichm\"uller
space.  We may cut $S$ into $2g-2$ pairs of pants so that the
boundary of each pair of pants is an $\rp^2$ geodesic, i.e.\ is a
straight line in an projective coordinate chart. Around each of
these geodesic loops is a holonomy action, which may be
represented as a matrix $H$ in $\Sl{3}$ once we choose an
appropriate frame.  The conjugacy class of this matrix is
invariant of our choice and constitutes the analog of
Fenchel-Nielsen's length parameter.  The eigenvalues of the
holonomy matrix $H$ are real, positive, and distinct.  So the
holonomy type may be described by the set of eigenvalues
$$\{\alpha_1,\alpha_2,\alpha_3\},\quad \alpha_i\in\re^+,\quad
\alpha_1\alpha_2\alpha_3=1,\quad \alpha_i\neq\alpha_j\mbox{ for
}i\neq j.
$$ Matrices in $\Sl{3}$ whose eigenvalues satisfy these
conditions are called \emph{hyperbolic}.  Order the eigenvalues
$\alpha_i$ so that the standard form of a hyperbolic holonomy
matrix is
\begin{equation} \label{standard-hyp-matrix}
\mathbf{D}(\alpha_1,\alpha_2,\alpha_3)\quad
\alpha_1>\alpha_2>\alpha_3>0
\end{equation}
(Here $\mathbf{D}$ denotes the diagonal matrix.)

\begin{figure}
\begin{center}
\scalebox{.3}{\includegraphics{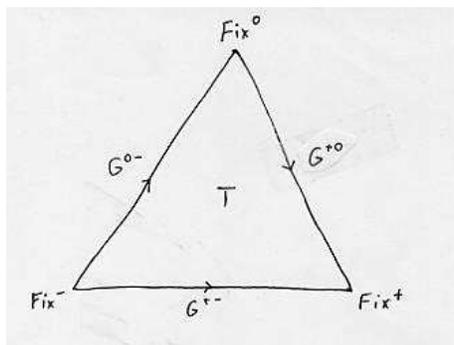}}
\end{center}
\caption{The Principal Triangle \label{principal-triangle}}
\end{figure}

There are three fix points of this hyperbolic holonomy matrix:
Fix$^+=[1,0,0]$ is attracting, Fix$^-=[0,0,1]$ is repelling, and
Fix$^0=[0,1,0]$ is a saddle fixed point.  The holonomy acts on
each coordinate line in $\rp^2$, and these lines split $\rp^2$
into four triangles.  Fix one open triangle $T$ to be the
projection of the first octant in $\re^3$ to $\rp^2$.  The
boundary of $T$ is formed by three line segments: $G^{+0}$
connects Fix$^+$ to Fix$^0$, $G^{0-}$ connects Fix$^0$ to Fix$^-$,
and $G^{+-}$ connected Fix$^+$ to Fix$^-$.  $G^{+-}$ is called the
\emph{principal} geodesic segment for the holonomy action.  See
Figure \ref{principal-triangle}.

We will also need to consider two other types of holonomy which
are degenerate in that they cannot occur in closed $\rp^2$
surfaces. The first is \emph{quasi-hyperbolic} holonomy, in which
the holonomy matrix is conjugate to
$$
\left( \begin{array}{ccc}
\alpha_1 & 1 & 0 \\ 0 & \alpha_1 & 0 \\ 0 & 0 & \alpha_3
\end{array} \right)\quad
\mbox{for}\quad \alpha_i>0,\quad \alpha_1^2\alpha_3=1, \quad
\alpha_1\neq\alpha_3.
$$ Also we will consider \emph{parabolic}
holonomy, in which the holonomy matrix is conjugate to
$$\left(
\begin{array}{ccc}
1&1&0 \\ 0&1&1 \\ 0&0&1
\end{array}
\right).$$ Choi \cite{choi94b} discusses all the types of holonomy
which can appear in an oriented $\rp^2$ surface.  We describe the
dynamics of quasi-hyperbolic and parabolic holonomies is
Subsections \ref{quasi-h-subsec} and \ref{parabolic-subsec}
respectively.

There are also twist parameters analogous to those in
Fenchel-Nielsen theory.  We recall the discussion in Goldman
\cite{goldman90a}. Around any oriented simple homotopically
nontrivial loop in a compact, convex $\rp^2$ surface $S=S_0$, the
holonomy type is hyperbolic. Inside the homotopy class of such a
loop, there is a unique representative which is a simple closed
principal geodesic loop $\mathcal L$. Then we may cut the surface
$S_0$ along $\mathcal L$ to form a possibly disconnected $\rp^2$
surface $S_0^{\rm cut}$ with principal geodesic boundary.  For
simplicity, we discuss only the case where $S_0^{\rm cut}= S^a_0
\sqcup S^b_0$ is disconnected. The other case is similar. Choose
$\rp^2$ coordinates on $S_0^a$ near the principal boundary
geodesic so that the holonomy matrix $H$ is in the standard form
(\ref{standard-hyp-matrix}), the principal boundary geodesic
$\mathcal{L}$ develops to the standard $G^{+-}$, and the developed
image of $S_0^a$ does not intersect the interior of $T$. Then the
other component $S_0^b$ is attached along $G^{+-}$ by placing the
image of its developing map inside $T$ across $G^{+-}$ from ${\rm
dev} \widetilde{S_0^a}$. (The inverse image of $\mathcal L$ in
the universal cover $\Omega=\widetilde{S_0}$ consists of not
just the line segment $G^{+-}$, but also a line segment
$\rm{hol}(\beta) G^{+-}$ for each $\beta$ in the coset space
$\pi_1S_0 / \langle \gamma\rangle$, where $\gamma$ is the element in
$\pi_1S_0$ determined by $\mathcal L$.  We must do a similar gluing
across each of these line segments.  For simplicity, we focus on
just the gluing across $G^{+-}$.)

The $\rp^2$
structure on the glued surface $S_0=S^a_0\cup_{\mathcal L} S^b_0$ is
then determined by an orientation-reversing projective
diffeomorphism $J$ from a neighborhood $N^a \subset S^a_0$ of
$\mathcal L$ to a similar neighborhood $N^b \subset S^b_0$.  In
terms of coordinates in the developed image near $G^{+-}$, $J$ may
be represented as the diagonal matrix $\mathbf{D}(1,-1,1)$, which
commutes with the holonomy matrix $H$.

\begin{figure}
\begin{center}
\scalebox{.5}{\includegraphics{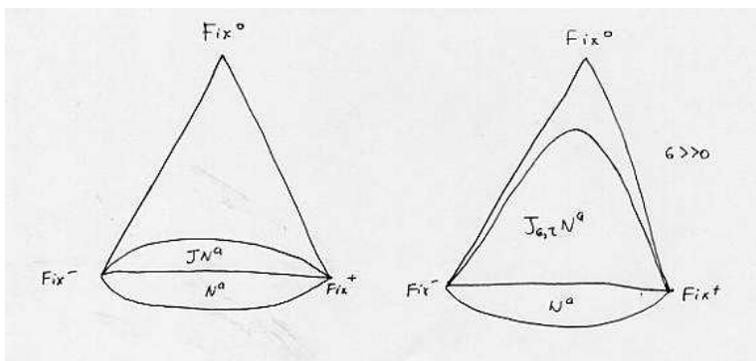}}
\end{center}
\caption{Twist Parameters \label{twist-triangle}}
\end{figure}

Now the twist parameters come in.  For real $(\sigma,\tau)$
consider the twist matrix
$$ \mathbf{M}(\sigma,\tau)= \mathbf{D}(e^{-\sigma-\tau},e^{2\tau},
e^{\sigma-\tau}).
$$ Then we form a new $\rp^2$ surface
$S_{\sigma,\tau}$ by gluing the neighborhood $S^b_0$ by the
projective involution $J_{\sigma,\tau}= \mathbf{M}(\sigma,\tau)J$
instead of the standard $J$. See Figure \ref{twist-triangle}.  Let
dev$_{\sigma,\tau}$ be equal to the standard developing map on
$\widetilde{S^a_0}$ as in the previous paragraph, and extend to
all of $\widetilde{S_{\sigma,\tau}}$ by using the gluing map
$J_{\sigma,\tau}$. We adapt Kim's terminology in \cite{kim99} to
call $\sigma$ the \emph{horizontal twist parameter} and $\tau$ the
\emph{vertical twist parameter}.  For the $\rp^2$ structure
determined by a hyperbolic Riemann surface, $\sigma$ corresponds
to the usual Fenchel-Nielsen twist parameter.

Note that as the vertical twist parameter $\tau\to+\infty$, the
image of the developing map
dev$_{\sigma,\tau}(\widetilde{S_{\sigma, \tau}})$ expands to
include all the interior of the principal triangle $T$. In this
case we have attached the entire \emph{principal half-annulus}
$\mathcal A =T/\langle H \rangle$ to $S^a_0$ along the principal
boundary geodesic $\mathcal L$.  There is only one way to attach
this principal half-annulus as an $\rp^2$ surface without
boundary. (Although there are two distinct ways to put a
non-principal geodesic boundary on $S^a_0\cup_{\mathcal L}\mathcal
A$---see e.g.\ Choi-Goldman \cite{choi-goldman97}.)

Similarly if $\tau\to-\infty$, dev$_{\sigma,\tau}
(\widetilde{S_{\sigma,\tau}})$ shrinks so that the glued part of
dev$_{\sigma,\tau}(\widetilde{S^b_0})$ vanishes, and
dev$_{\sigma,\tau} (\widetilde{S_{\sigma,\tau}})$ does not
intersect the principal triangle $T$. In this case, the principal
geodesic segment $\mathcal L$ is a natural boundary for the limit
$\rp^2$ surface.

Goldman \cite{goldman90a} also introduces interior parameters
associated to each pair of pants.  See also Kim \cite{kim99}. We
hope to use the methods of this paper to analyze them in future
work.

\section{Hyperbolic affine spheres and convex $\rp^n$ structures}
\label{h-a-s}

Recall the standard definition of $\rp^n$ as the set of lines
through 0 in $\re^{n+1}$.  There is a map
$P\!:\re^{n+1}\setminus0\to\rp^n$ with fiber $\re^*$.  For a
convex domain $\Omega\subset\re^n\subset\rp^n$ as above, then
$P^{-1}(\Omega)$ has two connected components.  Call one such
component $\mathcal{C}(\Omega)$, the \emph{cone over $\Omega$}.
Then any representation of a group $\Gamma$ into $\pgl{n+1}$ so
that $\Gamma$ acts discretely and properly discontinously on
$\Omega$ lifts to a representation into
$$\mathbf{SL}^\pm(n+1,\re)
= \{A\in\mathbf{GL}(n+1,\re):\det A=\pm1 \}$$
which acts on $\mathcal{C}(\Omega)$.
See e.g.\ \cite{loftin01}.

For a properly convex $\Omega$, then there is a unique hypersurface
asymptotic to the boundary of the cone $\mathcal{C}(\Omega)$ called
the hyperbolic affine sphere \cite{calabi72,cheng-yau77,cheng-yau86}.
This hyperbolic affine sphere $H\subset\mathcal{C}$ is invariant under
automorphisms of $\mathcal{C}(\Omega)$ in $\mathbf{SL}^\pm(n+1,\re)$.
The projection map $P$ induces a diffeomorphism of $H$ onto $\Omega$.
Affine differential geometry provides
$\mathbf{SL}^\pm(n+1,\re)$-invariant structure on $H$ which then
descends to $M=\Omega/\Gamma$.  In particular,
both the affine metric, which is a Riemannian metric conformal to
the (Euclidean) second fundamental form of $H$,
and a projectively flat connection whose geodesics are the $\rp^n$
geodesics on $M$, descend to $M$.
See \cite{loftin01} for details.  A fundamental fact
about hyperbolic affine spheres is due to Cheng-Yau \cite{cheng-yau86}
and Calabi-Nirenberg (unpublished):

\begin{thm} \label{cheng-yau-thm}
If the affine metric on a hyperbolic affine sphere $H$ is
complete, then $H$ is properly embedded in $\re^{n+1}$ and is
asymptotic to a convex cone $\mathcal C\subset \re^{n+1}$ which
contains no line.  By a volume-preserving affine change of
coordinates in $\re^{n+1}$, we may assume
$\mathcal{C}=\mathcal{C}(\Omega)$ for some properly convex domain
$\Omega$ in $\rp^n$.
\end{thm}

Below in Section \ref{wang-corr}, we recall a theory due to C.P.\ Wang
in the case $n=2$ and $M$ is oriented.  In this case, a properly
convex $\rp^2$ structure is given by certain data on a Riemann surface
$\Sigma$, and the developing map is given explicitly in terms of these
data by the solution to a first-order linear system of PDEs.

\section{Wang's developing map} \label{wang-corr}

C.P.\ Wang formulates the condition for a two-dimensional surface
to be an affine sphere in terms of the conformal geometry given by
the affine metric \cite{wang91}.  Since we rely heavily on this
work, we give a version of the arguments here for the reader's
convenience.  For basic background on affine differential
geometry, see Calabi \cite{calabi72}, Cheng-Yau \cite{cheng-yau86}
and Nomizu-Sasaki \cite{nomizu-sasaki}.

Choose a local conformal coordinate $z=x+iy$ on the hypersurface.
Then the affine metric is given by $h=e^{\psi}|dz|^2$ for some
function $\psi$. Parametrize the surface by $f:\mathcal{D}
\rightarrow \mathbb{R}^3$, with $\mathcal{D}$ a domain in
$\mathbb{C}$. Since $\{e^{-\frac{1}{2}\psi}
f_x,e^{-\frac{1}{2}\psi} f_y\}$ is an orthonormal basis for the
tangent space, the affine normal $\xi$ must satisfy this volume
condition (see e.g.\ \cite{nomizu-sasaki})
$$\det(e^{-\frac{1}{2}\psi} f_x,e^{-\frac{1}{2}\psi} f_y, \xi) = 1,$$
which implies
\begin{equation}
\det(f_z,f_{\bar{z}},\xi)=\sfrac{1}{2}i e^\psi. \label{det-eq}
\end{equation}

Now only consider hyperbolic affine spheres. By scaling in
$\mathbb{R}^3$, we need only consider spheres with affine mean
curvature $-1$.  In this case, we have the following structure
equations:
\begin{equation}
\left\{ \begin{array}{c}
D_X Y = \na_X Y + h(X,Y)\xi\\
D_X \xi = X
\end{array} \right. \label{struc}
\end{equation}
Here $D$ is the canonical flat connection on $\re^3$, $\na$ is a
projectively flat connection, and $h$ is the affine metric.
If the center of the affine sphere is $0$, then we also have $\xi=f$.

It is convenient to work with complexified tangent vectors, and we
extend $\na$, $h$ and $D$ by complex linearity.  Consider the
frame for the tangent bundle to the surface $\{ e_1 = f_z =
f_*(\frac{\partial}{\partial z}), e_{\bar 1}=f_{\bar z} =
f_*(\frac{\partial}{\partial {\bar z}}) \}$. Then we have
\begin{equation}
h(f_z,f_z)=h(f_{\bar z}, f_{\bar z})=0, \quad h(f_z,f_{\bar z}) =
\sfrac{1}{2}e^\psi.
\label{h-met}
\end{equation}
Consider $\theta$ the matrix of connection one-forms
$$\na e_i = \theta^j_i e_j, \quad i,j \in \{1,{\bar 1}\},$$
and $ {\hat \theta} $ the matrix of connection one-forms for the
Levi-Civita connection.   By (\ref{h-met})
\begin{equation}
{\hat \theta}^1_{\bar 1} ={\hat \theta}^{\bar 1}_1 = 0, \quad
{\hat \theta}^1_1 = \partial \psi, \quad
{\hat \theta}^{\bar 1}_{\bar 1} = {\bar \partial} \psi. \label{levi-cit}
\end{equation}

The difference ${\hat \theta} - \theta$ is given by the Pick form.  We
have
$${\hat \theta}^j_i - \theta^j_i = C^j_{ik} \rho^k,$$
where $\{ \rho^1 = dz,\rho^{\bar 1} = d{\bar z} \}$ is the dual
frame of one-forms.   Now we differentiate (\ref{det-eq}) and use
the structure equations (\ref{struc}) to conclude
$$\theta^1_1 + \theta^{\bar 1}_{\bar 1} = d \psi.$$
This implies, together with (\ref{levi-cit}), the apolarity condition
$$C^1_{1k} + C^{\bar 1}_{{\bar 1} k} = 0, \quad k \in \{1,{\bar 1} \}.$$
Then, when we lower the indices, the expression for the metric
(\ref{h-met}) implies that
$$C_{{\bar 1}1k} + C_{1{\bar 1}k} = 0.$$
Now $C_{ijk}$ is totally symmetric on three indices
\cite{cheng-yau86,nomizu-sasaki}. Therefore, the previous equation
implies that all the components of $C$ must vanish except
$C_{111}$ and $C_{ {\bar 1}{\bar 1}{\bar 1}} =
\overline{C_{111}}$.

This discussion completely determines $\theta$:
\begin{equation}
\left( \begin{array}{cc}  \theta^1_1 & \theta^1_{\bar 1} \\[1mm]
                                \theta^{\bar 1}_1 & \theta^{\bar 1}_{\bar1}
                  \end{array} \right)
 = \left( \begin{array}{cc}  \partial \psi & C^1_{{\bar 1}{\bar 1}}
                                d{\bar z} \\[1mm]
                                C^{\bar 1}_{11} dz & \bar{\partial} \psi
                  \end{array} \right)
=\left( \begin{array}{cc}  \partial \psi & \bar{U} e^{-\psi} d{\bar z} \\
                        U e^{-\psi} dz & \bar{\partial} \psi
                  \end{array} \right),
\label{conn-eq}
\end{equation}
where we define $U = C^{\bar 1}_{11} e^\psi$.

Recall that  $D$
is the canonical flat connection induced from ${\mathbb R}^3$.  (Thus,
for example, $D_{f_z}f_z = D_{\frac{\partial}{\partial z}} f_z =
f_{zz}$.)
Using this statement, together with (\ref{h-met}) and (\ref{conn-eq}),
the structure equations (\ref{struc}) become
\begin{equation}
\left\{ \begin{array}{c}
f_{zz} = \psi_z f_z + U e^{-\psi} f_{\bar z} \\
f_{{\bar z}{\bar z}} = {\bar U} e^{-\psi} f_z + \psi_{\bar z} f_{\bar z}
\\
f_{z{\bar z}} = \frac{1}{2}e^\psi f \end{array} \right. \label{fzz-eq}
\end{equation}
Then, together with the equations $f_z = f_z$ and $f_{\bar z} = f_{\bar
z}$, these form a linear first-order system of PDEs in $f$, $f_z$ and
$f_{\bar z}$:
\begin{eqnarray} \label{z-deriv}
\frac{\partial}{\partial z}
\left( \begin{array}{c} f \\ f_z \\ f_{\bar z} \end{array} \right)
&=&
\left( \begin{array}{ccc}
0&1&0 \\
0& \psi_z & U e^{-\psi} \\
\frac12 e^\psi & 0 & 0
\end{array} \right)
\left( \begin{array}{c} f \\ f_z \\ f_{\bar z} \end{array} \right),
\\ \label{zbar-deriv}
\frac{\partial}{\partial \bar z}
\left( \begin{array}{c} f \\ f_z \\ f_{\bar z} \end{array} \right)
&=&
\left( \begin{array}{ccc}
0&0&1 \\
\frac12 e^\psi & 0 & 0 \\
0 & \bar{U} e^{-\psi} & \psi_{\bar z} \\
\end{array} \right)
\left( \begin{array}{c} f \\ f_z \\ f_{\bar z} \end{array} \right).
\end{eqnarray}
In order to have a solution of the system (\ref{fzz-eq}),
the only condition is
that the mixed partials must commute (by the Frobenius theorem). Thus we
require
\begin{eqnarray}
\psi_{z {\bar z}} + |U|^2 e^{-2\psi} - \sfrac{1}{2} e^\psi  &=& 0,
\label{psi-eq} \\
\nonumber
U_{\bar z} &=& 0.
\end{eqnarray}

The system (\ref{fzz-eq}) is an initial-value problem, in that
given (A) a base point $z_0$, (B) initial values $f(z_0)\in\re^3$,
$f_z(z_0)$ and $f_{\bar z}(z_0)=\overline{f_z(z_0)}$, and (C) $U$
holomorphic and $\psi$ which satisfy (\ref{psi-eq}), we have a
unique solution $f$ of (\ref{fzz-eq}) as long as the domain of
definition $\mathcal{D}$ is simply connected.  We then have that
the immersion $f$ satisfies the structure equations (\ref{struc}).
In order for $f$ to be the affine normal of $f(\mathcal{D})$, we
must also have the volume condition (\ref{det-eq}), i.e.\ $\det
(f_z, f_{\bar z}, f)=\frac{1}{2}ie^\psi$.  We require this at the
base point $z_0$ of course:
\begin{equation}
\det(f_z(z_0), f_{\bar z}(z_0), f(z_0)) = \sfrac{1}{2}ie^{\psi(z_0)}.
\label{init-val}
\end{equation}
Then use (\ref{fzz-eq}) to show that the derivatives with respect
to $z$ and ${\bar z}$ of $\det(f_z,f_{\bar z},f)e^{-\psi}$ must
vanish. Therefore the volume condition is satisfied everywhere,
and $f(\mathcal{D})$ is a hyperbolic affine sphere with affine
mean curvature $-1$ and center $0$.

Using (\ref{fzz-eq}), we compute
$\det(f_z,f_{zz},f)=\frac{1}{2}iU$, which implies that $U$
transforms as a section of $K^3$, and $U_{\bar z}=0$ means it is
holomorphic.

Also, consider two embeddings $f$ and $\hat f$ from a
simply connected $\mathcal{D}$ to $\mathbb{R}^3$ which
satisfy (\ref{fzz-eq}) and the initial value condition (\ref{init-val})
for
some $z_0$ and $\tilde{z}_0$. Then consider the map $A \in
\gl{3}$
$$ A=\left( \begin{array}{c}
f(z_0) \\ f_z(z_0) \\ f_{\bar{z}}(z_0)
\end{array} \right)^{-1}
\left( \begin{array}{c}
\hat{f}(z_0) \\ \hat{f}_z(z_0) \\ \hat{f}_{\bar{z}}(z_0)
\end{array} \right)$$
  By the
volume condition (\ref{det-eq}), $A \in \Sl{3}$.  The uniqueness of
solutions to (\ref{fzz-eq}) then shows that $fA=\hat{f}$ everywhere.

We record all this discussion in the following
\begin{prop}[Wang \cite{wang91}]
\label{two-sphere} Let $\mathcal{D} \subset {\mathbb C}$ be a
simply connected domain. Given $U$ a holomorphic section of $K^3$
over $\mathcal{D}$, $\psi$ a real-valued function on $\mathcal{D}$
so that $U$ and $\psi$ satisfy (\ref{psi-eq}), and initial values
for $f$, $f_z$, $f_{\bar z}$ which satisfy (\ref{init-val}), we
can solve (\ref{fzz-eq}) so that $f(\mathcal{D})$ is a hyperbolic
affine sphere of affine mean curvature $-1$ and center $0$.  Any
two such $f$ which satisfy (\ref{fzz-eq}) are related by a motion
of $\Sl{3}$.
\end{prop}

More generally, if $\Sigma$ is
a Riemann surface with metric $h = e^{\phi} |dz|^2$.  Now write the
affine metric as
$e^{\psi} |dz|^2=e^u\,h$. Therefore,
$u$ is a globally defined function on $\Sigma$ and locally $\psi = \phi
+u$.  The Laplacian
$\Delta = 4 e^{-\phi} \partial_z \partial_{\bar{z}}$.  Therefore, $\psi$
solves (\ref{psi-eq})
exactly if the following equation in $u$ holds:
\begin{eqnarray}
\Delta u &=& 4 e^{-\phi}\psi_{z {\bar{z}}} - \Delta\phi  \nonumber\\
&=& 2 e^{-\phi} ( -2 e^{- 2 \psi} |U|^2 + e^{\psi} ) + 2\kappa  \nonumber
\\
&=&-4 e^{-2u} \|U\|^2 + 2 e^u + 2\kappa \label{wang-eq}
\end{eqnarray}
Here $\| \cdot \|^2 = | \cdot |^2 e^{-3\phi}$ denotes the metric on $K^3$
induced by $h$
and $\kappa=-\frac12 \Delta \phi$ is the curvature of $h$.

Note that this discussion gives an explicit description of the
developing map. Consider a Riemann surface $\Sigma$ equipped with
a holomorphic cubic differential $U$ and a conformal metric $h$.
Let $\mathcal{D}\subset\co$ be the universal cover of $\Sigma$. If
on $\Sigma$ there is  a solution $u$ to (\ref{wang-eq}) so that
$e^uh$ is complete, then Cheng-Yau and Calabi-Nirenberg's Theorem
\ref{cheng-yau-thm} above implies that the affine sphere
$f(\mathcal{D})$ is asymptotic to a convex cone $\mathcal{C}$
which contains no lines.  Therefore,
$\mathcal{C}=\mathcal{C}(\Omega)$ for a properly convex
$\Omega\subset \rp^2$.  As in Section \ref{h-a-s} above, the
projection map $P$ takes $f(\mathcal{D})$ diffeomorphically to
$\Omega$.  The developing map from $\mathcal{D}$ to $\Omega$ is
then explicitly $P(f)$, where $f$ satisfies the initial value
problem (\ref{fzz-eq}), (\ref{init-val}).

Consider as above $(\Sigma,U,e^uh)$ with $e^uh$ complete, and
$\mathcal{D}$ the universal cover of $\Sigma$.  For
$z_0\in\mathcal{D}$, choose a particular solution
$f\!:\mathcal{D}\to \re^3$ to the initial value problem
(\ref{fzz-eq}), (\ref{init-val}). Let $\gamma\in\pi_1M$ be a deck
transformation of $\mathcal{D}\to\Sigma$, which we take to be a
holomorphic automorphism of $\mathcal{D}$.  Then the uniqueness of
the hyperbolic affine sphere and of the initial value problem
imply that $f(\mathcal{D})=\gamma^*f(\mathcal{D})$.  Moreover, the
complexified frame in $\re^3$ $\{f,f_z,f_{\bar z}\}$ pulls back
under $\gamma$ to
\begin{equation} \label{gamma-frame}
\{f\circ\gamma,\gamma'f_z\circ\gamma,
\overline{\gamma'}f_{\bar z}\circ\gamma\}.
\end{equation}
In the particular cases considered below in Sections \ref{holo-sec}
and \ref{detail-sec}, the deck transformation $\gamma$ is of the
form
$$\gamma\!:z\mapsto z+c
$$
for a constant $c$.  Thus by (\ref{gamma-frame}), $\{f,f_z,f_{\bar
z}\}$ makes sense as a frame of a natural vector bundle over
$\mathcal{D}/\langle\gamma\rangle$.  Define the matrix $H_\gamma$
by
$$H_\gamma\!: \{f(z_0),f_z(z_0),f_{\bar z}(z_0)\} \mapsto
\{f(\gamma(z_0)), f_z(\gamma(z_0)), f_{\bar z}(\gamma(z_0)) \}.
$$
$H_\gamma$ maps the affine sphere $f(\mathcal{D})$ to itself and
satisfies $\det H_\gamma=1$.  Then $H_\gamma$ is conjugate to a
matrix in $\Sl{3}$---simply use the real frame $\{f,f_x,f_y\}$
instead---and, by projecting to $\pgl{3}$, determines the holonomy
of the $\rp^2$ structure along a loop in $\Sigma$ whose endpoints
lift to $z$ and $\gamma(z)$. We record this in

\begin{prop} \label{calc-holo-prop}
Consider $(\Sigma,U,e^uh)$ as above so that $e^uh$ is complete,
and let $\mathcal{D}\subset\co$ be the
universal cover of $\Sigma$.  If a loop in $\Sigma$ can be represented
by a deck transformation of the form $\gamma(z)=z+c$, then the frame
$\{f,f_z,f_{\bar z}\}$ may be used to calculate the holonomy around this
loop.  A matrix in the conjugacy class of the
holonomy may be obtained by integrating the
initial-value problem (\ref{z-deriv}-\ref{zbar-deriv}),
(\ref{init-val}) along a path whose
endpoints in $\mathcal{D}$ are $z$ and $z+c$.
\end{prop}

We remark that in the particular case the initial metric
$h$ is hyperbolic (i.e.\ with constant
curvature $-1$), we have the equation
$$
\Delta u = -4 e^{-2u} \|U\|^2+2e^u-2,
$$
which has a unique solution on a compact Riemann surface
$\Sigma$ of genus $g\ge2$ for any
$U\in H^0(\Sigma,K^3)$ (Wang, \cite{wang91}).  In the next section,
we extend this result to noncompact Riemann surfaces which admit a
hyperbolic metric of finite volume.

\section{Finding solutions} \label{finding-sol}

\subsection{The Ansatz}
Consider $\Sigma=\bar{\Sigma}\setminus \{p_i\}$ be a Riemann surface of
finite type equipped with a complete hyperbolic metric.  Consider
$U$ a section of $K_{\Sigma}^3$ with poles of order at most three allowed
at the punctures $p_i$.  In other words $U \in
H^0(K_{\bar{\Sigma}}^3\otimes
\prod_i [p_i]^3)$.
We want to find a metric $h$ so that
\begin{equation}
-4\|U\|^2+ 2+ 2\kappa \to 0\quad \mbox{at the} \quad p_i,
\label{asymp-eq}
\end{equation}
 so that $u=0$ is an approximate solution to (\ref{wang-eq}).

Near a puncture point $p_i$, consider $z=z_i$ the local coordinate
function so that the hyperbolic metric is exactly
\begin{equation}
h=\frac4{|z|^2(\log|z|^2)^2}\,|dz|^2 \label{hyp-metric}
\end{equation}
near the puncture $\{z=0\}$. (For now we drop the notational
dependence on $i$.)  Such a coordinate $z$ is called a \emph{cusp
coordinate} near the puncture. Cusp coordinates are unique up to a
rotation $\tilde z= e^{i\theta}z$. Near the puncture,
$U=R\,z^{-3}dz^3+O(z^{-2})$ for a complex number $R$. We call $R$
the residue of $U$ at the puncture.  If $R=0$, then we just leave
the hyperbolic metric, and
$\|U\|^2=O\left(|z|^2(\log|z|^2)^6\right)$; therefore,
(\ref{asymp-eq}) is satisfied.

For $R\neq 0$, however, we choose a flat metric near the puncture.  Let
\begin{equation}
h=\frac{2^{\frac13}|R|^{\frac23}}{|z|^2}\, |dz|^2 \label{flat-metric}
\end{equation}
near $z=0$.  This metric then satisfies the asymptotic requirement
(\ref{asymp-eq}).

Now for a given $U$, we can patch these metrics together on $\Sigma$ by
requiring that $h$ be hyperbolic outside of a neighborhood of those $p_i$
for which $R_i\neq 0$.  In particular, $h$ must be hyperbolic on a
neighborhood of all the zeros of $U$. In a neighborhood of each $p_i$
for which $R_i\neq0$, we make $h$ be the flat metric (\ref{flat-metric}).
On the remainder of $\Sigma$ (which consists of annular necks around
each $p_i$ with nonzero residue), we let $h$ be an arbitrary conformal
metric smoothly interpolating the flat and hyperbolic metrics.  Note that
by this construction, we have two types of punctures $p_i$.
We name these ends according to the holonomy of the $\rp^2$ surface
we will construct from $(\Sigma,U,h)$---see Table \ref{holonomy-table} below.
We call those punctures $p_i$ for which
$R_i\neq0$ the QH ends of $(\Sigma,h)$, since
the $\rp^2$ holonomy will be quasi-hyperbolic or hyperbolic according to
whether Re$\,R_i=0$ or not.
Those $p_i$ for which
$R_i=0$ are the parabolic ends of $(\Sigma,h)$.

Here is a more explicit description of the metric near a
QH end.
In the conformal coordinate $z_i$ as above, we define
\begin{equation} \label{def-h}
h=\left\{ \begin{array}{c@{\mbox{ for }}c}
\frac{2^{\frac13}|R|^{\frac23}}{|z|^2}\, |dz|^2 & |z_i|<c_i \\
e^{\rho_i} |dz|^2 & c_i\le|z_i|\le C_i \\
\frac4{|z|^2(\log|z|^2)^2}\,|dz|^2 & |z_i|>C_i
\end{array} \right.
\end{equation}
Here $c_i<C_i$ are appropriate radii and $e^{\rho_i}$ is a smooth
interpolation between the two metrics. We require the zeros of $U$ to
be away from the QH
ends of $(\Sigma,h)$, so that for $0<|z_i|\le C_i$,
$\|U(z_i)\|\ge \delta_i>0$.  This is possible since at each
QH end, $\lim_{z_i\to0}\|U(z_i)\|= \frac1{\sqrt 2}$.

\subsection{Solving Wang's equation} \label{find-sol}
Now we will find solutions to (\ref{wang-eq}) for the given $\Sigma$ and
$U$, and the metric $h$ constructed in the previous section.  We will
construct barriers on $\Sigma$ to show that the solution $u$ we find will
be bounded and will approach zero near the QH
ends of $(\Sigma,h)$.

To find a supersolution to (\ref{wang-eq}), define
\begin{equation}
L(u)=\Delta u +4e^{-2u}\|U\|^2-2e^u-2\kappa \label{operator}
\end{equation}
so that $L(u)=0$ is our equation.  Then near a QH
end of $(\Sigma,h)$,
consider $v=\beta|z|^{2\alpha}$ for $\alpha$, $\beta$ positive constants.
Calculate
\begin{eqnarray}\nonumber
L(v)&=& 4\beta\alpha^2(2^{-\frac13})|R|^{-\frac23}\,|z|^{2\alpha} +
[2+O(|z|)]e^{-2\beta|z|^{2\alpha}}-2e^{\beta|z|^{2\alpha}}\\
\label{L-eq}
&=&\big[4\alpha^2(2^{-\frac13})|R|^{-\frac23}-3\big]v + 2e^{-2v}-2e^v+3v+
e^{-2v}O(|z|).
\end{eqnarray}
If we choose $\alpha$ small enough, the first term is negative, and
we can check that $2e^{-2v}-2e^v+3v$ is negative for all $v>0$. The
term $e^{-2v}O(|z|)$ is dominated by the first term for $\alpha$ small,
$\beta$ large and $z$ near $0$. So $L(v)\le0$ on a neighborhood
$\mathcal{N}$ of $z=0$, and $\mathcal{N}$ can be made independent of
the choice of $\beta$ for $\beta\gg0$.

So consider a smooth positive function $f$ on $\Sigma$ which
satisfies $f=|z_i|^{2\alpha_i}$ on the neighborhood
$\mathcal{N}_i$ corresponding to each QH end of $(\Sigma,h)$ for a
suitably small $\alpha_i$.  Near the parabolic ends of $\Sigma$,
let $f$ be a positive constant,and let $f$ be smooth and positive
on all $\Sigma$. Then for $\beta$ large, $L(\beta f)\le0$ on all
of $\Sigma$, since the $-2e^u$ term in (\ref{operator}) dominates
outside the $\mathcal{N}_i$.  This will be our supersolution
$S=\beta f$.  Note that $S$ is bounded and positive, and $S\to0$
at each QH end of $(\Sigma,h)$.

Finding a subsolution is somewhat more delicate, since the
presence of zeros of $U$ means that the positive term
$4e^{-2u}\|U\|^2$ in (\ref{operator}) cannot dominate all the
others for $u \ll 0$.  We will look to the curvature term
$-2\kappa$ instead for positivity. In particular, we have required
the metric to be hyperbolic (so $-2\kappa=2$) wherever $\|U\|$ is
small.

First, near each QH end, we can consider $w=-\beta |z|^{2\alpha}$
for $\alpha,\, \beta>0$. As above
$$ L(w)=\big[4\alpha^2(2^{-\frac13})|R|^{-\frac23}-3\big]w +
2e^{-2w}-2e^w+3w+e^{-2w}O(|z|),$$ and for negative $w$,
$2e^{-2w}-2e^w+3w>0$. For small $\alpha$, large $\beta$, the first
term is positive and dominates the term $e^{-2w}O(|z|)$.
Therefore, as above, we have neighborhoods of the QH ends
$\mathcal{N}_i$ which do not depend on $\beta$ for $\beta\gg0$,
and $L(w)\ge 0$ on these $\mathcal{N}_i$.

Now recall the situation in equation (\ref{def-h}).  Near a QH
end, we have the metric is hyperbolic for $|z_i|>C_i$, flat for
$|z_i|<c_i$. Also, we assume that the $\mathcal{N}_i \subset
\{|z_i|<c_i\}$. We have no control over the curvature $\kappa$ for
$c_i\le|z_i|\le C_i$, but we do know that $\|U\|\ge \delta_i$
there.  Therefore, we can let $\beta_i$ become large so that the
term $4e^{-2w}\|U\|^2$ dominates the others on $\{|z_i|\le
C_i\}\setminus \mathcal{N}_i$, and thus
$L(-\beta_i|z_i|^{2\alpha_i})\ge0$ for $|z_i|\le C_i$.

By making some $\beta_i$ larger if necessary, we can make sure that the
values of $-\beta_i|z_i|^{2\alpha_i}$ are all equal to some negative
constant $-B$ on the circles
$\{|z_i|=C_i\}$.  Then we define the subsolution $s$ as
$$ s=\left\{
\begin{array}{c@{\quad}c}
-\beta_i|z_i|^{2\alpha_i} & \mbox{on each } |z_i|\le C_i \\
-B & \mbox{elsewhere}
\end{array}
\right.$$ Then in the hyperbolic part of $(\Sigma,h)$, $L(s)=
4e^{2B}\|U\|^2-2e^{-B}+2>0$.  On the circles $\{|z_i|=C_i\}$, $s$
is not smooth, but since $\Delta(s)\ge0$ as a distribution there,
$s$ is a suitable lower barrier.  So $L(s)\ge0$ on $\Sigma$ and
$s\to0$ at the QH ends of $(\Sigma,h)$.  Also note $s$ is bounded
and negative.

Now that we have upper and lower barriers, we can find a solution
to (\ref{wang-eq}) on $(\Sigma,h)$.  Write $\Sigma = \bigcup_j
\Omega_j$, where the $\Omega_j$ are a sequence of compact
submanifolds with boundary which exhaust $\Sigma$. Then on each
$\Omega_j$, we can solve the Dirichlet problem $L(u)=0$ on
$\Omega_j$ and $u=0$ on $\partial \Omega_j$ (as in e.g.\
\cite{gilbarg-trudinger}, Thm.\ 17.17; the main thing to check
here is that the nonlinear operator $L$ is decreasing as a
function of $u$.)  Call this solution $u_j$.  By the maximum
principle, we have $S\ge u_j \ge s$.

These bounds on the $u_j$ then give uniform local $L^p$ bounds on
the right hand side of (\ref{wang-eq}), and therefore by the
elliptic theory \cite{gilbarg-trudinger}, we have local $W^{2,p}$
bounds on the $u_j$. This is enough to ensure that a subsequence
of the $u_j$ converges uniformly to a solution $u$ on $\Sigma$.
Higher regularity of $u$ is standard, and the barriers $S$ and $s$
ensure that $S\ge u\ge s$. Thus $u\to0$ at the QH ends of
$(\Sigma,h)$ and $u$ is bounded everywhere.

We can also show, using Cheng and Yau's maximum principle for complete
manifolds \cite{cheng-yau75}, that the $u$ we have constructed is the
unique bounded solution to (\ref{wang-eq}).

\begin{prop} \label{unique}
There is only one bounded solution to (\ref{wang-eq}) for a given
$(\Sigma,U)$ and metric $h$ as constructed above.
\end{prop}

\begin{proof}
If $u$ and $\tilde u$ are two solutions to (\ref{wang-eq}) so that
$|u|,|\tilde u| \le M$, then $u-\tilde u$ satisfies
$$ \Delta (u-\tilde u) = g(x,u)-g(x,\tilde u),$$
where $g(x,u)=-4e^{-2u}\|U\|^2+2e^u$ is strictly increasing in $u$.
There is a positive constant
$$ C=\inf\{\partial_u g(x,u) : x\in \Sigma, \, u\in[-M,M] \}$$
so that
$$ \Delta (u-\tilde u) \ge C(u-\tilde u).$$
Since $(\Sigma, h)$ is complete and
has bounded Ricci curvature, Cheng and Yau's result implies
$$\forall\, \epsilon >0, \quad \exists\, x_\epsilon \in \Sigma \quad
\mbox{so that}\quad \Delta(u-\tilde u)(x_\epsilon)\le \epsilon,\quad
(u-\tilde u)(x_\epsilon) \ge B-\epsilon,$$
where $B=\sup_\Sigma (u-\tilde u)$.
Therefore,
$$ \epsilon \ge \Delta(u-\tilde u)(x_\epsilon) \ge C(u-\tilde u)(x_\epsilon)
\ge C(B-\epsilon)$$
Then $B \le \epsilon\frac{1+C}C$ for all $\epsilon >0$,
and thus $B\le0$.  A similar argument
shows $\inf_\Sigma(u-\tilde u)\ge0$ also. So $u=\tilde u$ on all of $\Sigma$.
\end{proof}

Finally we will need bounds on the gradient of $u$.  We use the
$L^p$ theory again to accomplish this.

\begin{lem} \label{gradient-bound}
Let $|\nabla\cdot|$ denote the norm of the gradient with respect to the
metric $h$ and let $f=-4 e^{-2u} \|U\|^2 + 2 e^u + 2\kappa$ be the right
hand side of (\ref{wang-eq}).  Then there is a constant $K$ independent
of $x\in \Sigma$ such that
$$  |\nabla u(x)|\le K (\|u\|_x + \|f\|_x),$$
where $\|\cdot\|_x$ denotes the sup norm in a geodesic ball of radius 1
around $x$.
\end{lem}

This lemma immediately shows that $|\nabla u|$ is always bounded and it
approaches zero at the QH
ends of $(\Sigma,h)$, since $u\to0$ at a QH
end implies $f\to0$ there as well by (\ref{flat-metric}).

\begin{proof}  Since the ends of $\Sigma$ are constant curvature $0$ or
$-1$, $\Sigma$ has bounded geometry.  In other words, there are uniform
constants $A<1$, $B_n$ so that for any $x\in\Sigma$,
\begin{itemize}
\item
There is a quasi-coordinate ball $\mathcal{B}$ of radius $A$ around $x$.
(Take some neighborhood $\mathcal{N}$ of $x$ in $\Sigma$, and pull back
the metric to the universal cover $\tilde{\mathcal{N}}$ of $\mathcal{N}$.
Our quasi-coordinate ball $\mathcal{B}\subset\tilde{\mathcal{N}}$
is a geodesic ball of radius $A$
centered at a lift of $x$ and
properly contained in $\tilde{\mathcal{N}}$.)
In these coordinates in $\mathcal{B}$, we have
\item
The metric $g_{ij}$ satisfies $|g_{ij}-\delta_{ij}|<B_0$.
\item
The ordinary $n^{\rm th}$ derivatives of $g_{ij}$ are bounded by $B_n$.
\end{itemize}
The usual geodesic normal coordinate balls satisfy these conditions.
These are the conditions we need to apply the $L^p$ estimates.

Choose $p>2$.
By the elliptic theory \cite{gilbarg-trudinger}, for uniform constants
$C$, $C'$, and a smaller ball $\mathcal{B}'$, also centered at $\tilde x$, we
have
$$ |\nabla u (x)|\le \|u\|_{C^1(\mathcal{B}')} \le
C\|u\|_{W^{2,p}(\mathcal{B}')} \le
C'\left(\|u\|_{L^p(\mathcal{B})}+ \|f\|_{L^p(\mathcal{B})}\right).$$
The second inequality follows by the Sobolev embedding theorem and the
third by interior $L^p$ estimates. The $L^p$ norm is in turn dominated by
the sup norm as required.
\end{proof}

We record the above discussion in a proposition.

\begin{prop} \label{bound-prop}
Given $\Sigma$, $U$ and $h$ as above, there is a unique bounded solution
$u$ to (\ref{wang-eq}).  $u$ is smooth and approaches zero at the
QH ends
of $(\Sigma,h)$.  Furthermore, the norm of the gradient $|\nabla u|$ is
bounded and approaches zero at the QH ends of $(\Sigma,h)$.
Specifically,
near each QH end of $(\Sigma,h)$, there are constants
$\alpha_i,\,\beta_i>0$
so that $|u|,\,|\nabla u| \le \beta_i|z_i|^{2\alpha_i}$. The metric
$e^uh$
is complete, and
$(\Sigma,e^uh,U)$ determine a convex $\rp^2$ structure on the surface.
\end{prop}

\begin{proof}
We have already proved all but the last sentence.  The affine
metric $e^uh$ is complete since $u$ is bounded and $h$ is
complete.  The statement about $\rp^2$ structures follows from
Wang's work on affine spheres as above, and Cheng and Yau's
classification of affine spheres with complete affine metric
\cite{cheng-yau86}.  See \cite{labourie97} or \cite{loftin01} for
more details about affine spheres and $\rp^n$ structures.
\end{proof}

Note that in \cite{labourie97,loftin01}, it is shown that a convex
$\rp^2$ structure on a compact oriented surface $S$ is equivalent
to a pair $(\Sigma,U)$ of $\Sigma$ a conformal structure on $S$
and $U$ a cubic differential on $\Sigma$.

\section{Holonomy type of the ends} \label{holo-sec}

In this section, we will use the asymptotics of the affine metric
$e^uh$ computed above and Wang's integrable system for the
associated affine sphere to compute the asymptotics of the $\rp^2$
structure on $\Sigma$ near each of its ends. The cases of the QH
and parabolic ends of $\Sigma$ will be treated separately.

\subsection{Topological setup} \label{setup}
Represent the universal cover of $\Sigma$ as the upper half-plane
$\mathbb{H}=\{w=x+iy\in \mathbb{C}:y>0\}$. Each puncture of
$\Sigma$ corresponds to a parabolic element of Aut$(\mathbb{H})$,
which in turn is conjugate to the map $\gamma\!:w\mapsto w+2\pi$.
For a given puncture $p$, we have covering maps
$$\mathbb{H}\stackrel{\zeta}{\to}\mathbb{D}_0\stackrel{\xi}{\to}\Sigma$$
Here $\mathbb{D}_0$ is the punctured disk $\{z\in \mathbb{C}:0<|z|<1\}$,
and $\xi$ extends to map $z=0$ to the puncture $p$.  Also, we define
$\zeta(w)=e^{iw}$; then the map $\gamma$ generates the deck
transformations
for the covering map $\zeta$.  Recall that on $\Sigma$
near the end $z=0$, $U= \frac{R}{z^3}[1+O(z)]dz^3$, and so
\begin{equation} \label{U-upstairs}
\zeta^*U=-iR[1+O(e^{-y})]dw^3.
\end{equation}

We consider neighborhoods of the puncture of the form $N =
\{0<|z|<\epsilon\}$, which is topologically a cylinder.  Below we
will consider explicit paths from a base point $p\in N$ to the
puncture.  As above, lift $N$ to the set $\tilde{N}=\{w:
y>-\log\epsilon\}$, and the base point $p$ to its lift in $\tilde
N$.  We will consider particular paths in $\tilde{N}$
corresponding to rays along which $y\to\infty$.  Pushed down to
the $z$ coordinate, these paths go to the puncture $z=0$.

Solving the initial value problem (\ref{fzz-eq}), as determined by
$(\Sigma,U,e^uh)$ then provides a developing map from the
universal cover $\mathbb H$ of $\Sigma$ into $\rp^2$.  Denote $S$
as the $\rp^2$ surface constructed in this way.  In $S$, $N$ is
topologically a cylindrical neighborhood of the end.  We consider
the $\pi_1(S)$ with respect to the basepoint $p\in N \subset S$.
We only consider paths from $p$ to the end which remain in $N$,
and thus we need only consider paths in the universal cover
$\tilde S$ which remain in dev$(\tilde{N})$. All this will supply
a very explicit model of the developing map near the end, upon
appropriate choice of coordinates on $\rp^2$.

\subsection{The main holonomy computation.} \label{main-holo-calc}
We present a simple argument to
calculate the holonomy around the puncture for these singular surfaces.

Consider the frame $\{f,f_w,f_{\bar w}\}$.  Then (\ref{fzz-eq})
shows that
\begin{equation} \label{x-hol}
\frac{\partial}{\partial x}
\left( \begin{array}{c}
f\\f_w\\f_{\bar w}
\end{array} \right)
=
\left( \begin{array}{ccc}
0&1&1\\
\frac12e^{\psi} & \psi_w& Ue^{-\psi} \\
\frac12e^{\psi} & \bar{U}e^{-\psi} & \psi_{\bar w}
\end{array} \right)
\left( \begin{array}{c}
f\\f_w\\f_{\bar w}
\end{array} \right)
\end{equation}
As above, $\psi=\phi+u$ and the initial metric $h=e^\phi|dw|^2$.  Define
${\bf A}_y$ to be the matrix in equation (\ref{x-hol}).

\begin{lem} \label{A-lemma}
$\D \lim_{y\to\infty}{\bf A}_y={\bf A}=
\left( \begin{array}{ccc}
0&1&1\\
2^{-\frac23}|R|^{\frac23} & 0 & -i 2^{-\frac13}R|R|^{-\frac23} \\
2^{-\frac23}|R|^{\frac23} & i2^{-\frac13}\bar{R}|R|^{-\frac23} & 0
\end{array} \right)$ \\
uniformly in $x$.  (If $R=0$, we put all the matrix entries involving
$R$ to be zero.) The characteristic polynomial of $\bf A$ is
\begin{equation} \label{char-poly}
\chi(\lambda)=\lambda^3-3(2^{-\frac23})|R|^{\frac23}\lambda - {\rm Im
}\,R.
\end{equation}
\end{lem}

\begin{proof}
In the case of a QH
end, $R\neq0$ and Proposition \ref{bound-prop}
shows
that in the $w$ coordinates, each of $u,u_w=O(e^{-2\alpha y})$ as
$y\to\infty$.  Then  the definition of $h$ (\ref{flat-metric}) and the
asymptotics of $U$ (\ref{U-upstairs}) provide the result.

If $R=0$, $h=\frac{|dw|^2}{y^2}$ and thus $\phi=-2\log y$. Proposition
\ref{bound-prop} then shows that $|u|\le C$ and $|u_w|\le\frac{C}y$.
These,
together with the asymptotics of $U$ (\ref{U-upstairs}), complete the
proof.
\end{proof}

We can immediately find the eigenvalues of any holonomy matrix
around each end. For a fixed $y\gg0$,  the loop $|z|=e^{-y}$ on
$\Sigma$ lifts to the line segment $(x,y)$ in $\mathbb{H}$, where
$x$ goes from $0$ to $2\pi$.  Let $H_y$ be the holonomy matrix of
the connection $D$ with respect to the frame $\{f,f_w,f_{\bar
w}\}$ around this loop. This is justified by Proposition
\ref{calc-holo-prop} above. Then $H_y=\Phi(2\pi)$, where $\Phi$
solves the initial value problem \{(\ref{x-hol}), $\Phi(0)=I$\}.
Since $D$ is flat and the loops are freely homotopic, all $H_y$
for $y\gg0$ are conjugate to each other in $\bf{GL}(3,\co)$.
Moreover, as the fundamental solution to (\ref{x-hol}),
\begin{equation} \label{hol-limit}
\lim_{y\to\infty}H_y\to e^{2\pi \bf{A}}.
\end{equation}
  This last statement
follows by Lemma \ref{A-lemma} and the theory of ODEs with
parameters \cite{hartman}. All the matrices $H_y$ for $y\gg0$ have
the same eigenvalues, and (\ref{hol-limit}) shows that they are
$e^{2\pi\lambda_i}$, where $\lambda_i$ are the eigenvalues of $\bf
A$.  Note that since Tr$\,{\bf A}=0$, $H_y$ is conjugate to a
matrix in $\Sl{3}$.

In the case of repeated roots of $\chi(\lambda)$, we {\em cannot}
conclude, however, that the $H_y$ have the same conjugacy type as
$e^{2\pi {\bf A}}$. As we'll see below, this is false, since the
matrix $e^{2\pi{\bf A}}$ in this case can be approximated by
matrices with the same eigenvalues but whose Jordan decomposition
consists of maximal Jordan blocks.  Indeed, we'll see below in
Subsection \ref{perturb} that the matrix through which $H_{y_0}$
and $H_y$ are conjugate diverges as $y\to\infty$ and $y_0$ is
fixed.

The discriminant of $\chi(\lambda)$ is
\begin{equation} \label{disc-eq}
 \mathcal{D} = -\sfrac14 |R|^2 + \sfrac14 (\rm{Im}\,R)^2.
\end{equation}
So $\mathcal{D}\le0$ always, and $\mathcal D =0$ only if Re$\,R=0$.
Therefore, $\chi(\lambda)$ only has real roots, and these are repeated
if and only if Re$\,R=0$.  So if Re$\,R\neq0$, we know the holonomy type
is given by the hyperbolic holomomy matrix whose eigenvalues are
$e^{2\pi\lambda_i}$ for $\lambda_i$ the eigenvalues $\bf A$.

\begin{prop} \label{hyperbolic-prop}
If Re$\,R\neq0$, then the holonomy around the puncture is conjugate
to
$$ \left( \begin{array}{ccc} e^{2\pi\lambda_1} & 0&0\\
0& e^{2\pi\lambda_2}&0 \\
0&0& e^{2\pi\lambda_3} \end{array} \right),$$
where $\lambda_i$ are the roots of $\chi(\lambda)$.  $\sum_i\lambda_i=0$,
and the $\lambda_i$ are real and distinct.
\end{prop}

\subsection{Quasi-hyperbolic holonomy.} \label{q-h-holo}
\begin{prop}\label{q-h-holo-prop}
Let Re$\,R=0$, but $R\neq0$, then the holonomy type is conjugate to
$$ \left( \begin{array}{ccc} e^{2\pi\lambda_1} & 1&0\\
0& e^{2\pi\lambda_1}&0 \\
0&0& e^{2\pi\lambda_3} \end{array} \right),$$
where $\lambda_i$ are the roots of $\chi(\lambda)$, with
$\lambda_1$ the repeated root.  $2\lambda_1+
\lambda_3=0$, $\lambda_1\neq\lambda_3$, and $\lambda_i\in\re$.
\end{prop}

\begin{proof} We have two choices for the holonomy:
$$
F= \left( \begin{array}{ccc} e^{2\pi\lambda_1} & 0 & 0 \\
 0 & e^{2\pi\lambda_1} & 0\\
 0 & 0 & e^{2\pi\lambda_3}
\end{array} \right)
\quad \mbox{or} \quad G= \left( \begin{array}{ccc}
e^{2\pi\lambda_1} & 1 & 0\\
0 & e^{2\pi\lambda_1} & 0 \\
0 & 0 & e^{2\pi\lambda_3}
\end{array} \right)
$$
Let $\alpha_i=e^{2\pi\lambda_i}$. A result of Choi \cite[Prop.\
2.3]{choi94b}, rules out the case of $F$ for a surface of negative
Euler characteristic. The result only applies, however, to a
surface whose end has the structure of an $\rp^2$ surface with
convex boundary. To find such a boundary, choose coordinates so
that $F$ is the lift of the element of the fundamental group
corresponding to holonomy around the end.  The developing image
$\Omega$ must contain a point $p=[x,1,z]$, written in homogeneous
coordinates in $\rp^2$. $\gamma^n$ then takes $p\mapsto
p_n=\left[x,1,\left(\frac{\alpha_3}{\alpha_1}\right)^n z \right]$.
All of these $p_n$ must be in $\Omega$, and by convexity, there
must be some line segment $\overline{p_n p_{n+1}}$ must be in
$\Omega$.  The action of powers of $\gamma$ ensure that all such
line segments are in $\Omega$; so the entire geodesic segment
$\{[x,1,sz] : s\in(0,\infty)\}$ is in $\Omega$.  On the quotient
surface $S=\Omega/\pi_1$, this is a geodesic loop isolating the
end from the rest of the surface. Cut along this geodesic loop and
then apply Choi's result to get a contradiction. Therefore, the
holonomy in this case is quasi-hyperbolic.
\end{proof}

\subsection{Parabolic holonomy.} Finally if $R=0$, then we have all the
eigenvalues of $e^{2\pi {\bf A}}$ are 1.

\begin{prop} \label{N-lemma}
Let $S$ be a properly convex $\rp^2$ surface.  In other words,
$S=\Omega/\Gamma$, where
$\Omega$ is a convex bounded open subset of some $\re^2 \subset \rp^2$,
and $\Gamma$ is a
subgroup of $\pgl{3}$ acting properly discontinuously on $\Omega$.  Any
element
$\gamma\in\Gamma$ whose set of eigenvalues is $\{1\}$ must be conjugate
to $N$, which consists of one $3\times3$ Jordan block.
\end{prop}

\begin{proof}
$\gamma$ must be conjugate to one of
$$I=
\left( \begin{array}{ccc}
1&0&0\\
0 &1&0 \\
0 &0&1
\end{array} \right),\quad
Q=
\left( \begin{array}{ccc}
1&1&0\\
0 &1&0 \\
0 &0&1
\end{array} \right), \quad
N=\left( \begin{array}{ccc}
1&1&0\\
0 &1&1 \\
0 &0&1
\end{array} \right) $$
The identity map $I$ is obviously not possible.  We now rule out $Q$.
Choose
coordinates so that $Q$ is the lift of $\gamma$ in $\Sl{3}$.  $\Omega$
must contain a point $p=[x,1,z]$, written in homogeneous coordinates on
$\rp^2$.  $\gamma^n$ then takes $p\mapsto p_n=[x+n,1,z]$, and so each
$p_n\in \Omega$.  Because $\Omega$ is convex in some $\re^2\subset\rp^2$,
it must then contain a line segment between two points $p_n$ and
$p_{n+1}$
for some $n$. By the action of powers of $\gamma$, it must contain all
such line segments.  In short, $\Omega$ contains the line
$\{[x+t,1,z]\! : t\in\re\}$.  As a {\em properly} convex domain
$\Omega$ cannot contain a whole line, $\gamma$ cannot be conjugate to
$Q$.
\end{proof}

\subsection{Results.} We record these results in Table
\ref{holonomy-table}.

\begin{table}

\caption{\label{holonomy-table} Residue and Holonomy}

\begin{tabular}{|c|l|l|}

\hline
 Residue & Holonomy type & Holonomy name\\
\hline
\hline
$R=0$ &
$\left(
\begin{array}{ccc}
1&1&0 \\ 0&1&1 \\ 0&0&1
\end{array}
\right)$
& Parabolic \\
\hline
\begin{tabular}{c}
Re$R=0$, \\
$R\neq0$
\end{tabular}
&
$\left(
\begin{array}{ccc}
\alpha_1 & 1 & 0 \\
0 & \alpha_2 & 0 \\
0&0&\alpha_3
\end{array}
\right) \quad (\alpha_2=\alpha_1)$
& Quasi-hyperbolic \\
\hline
Re$R\neq0$ &
$\left(
\begin{array}{ccc}
\alpha_1 & 0 & 0 \\
0 & \alpha_2 & 0 \\
0&0&\alpha_3
\end{array}
\right)$
& Hyperbolic \\
\hline
\end{tabular}

{Here $\chi\left((2\pi)^{-1}\log \alpha_i\right)=0$, $\prod\alpha_i=1$,
$\alpha_i>0$.}

\end{table}

Note that the holonomy type is uniquely determined by Im$\,R$ and $|R|$
alone,
by (\ref{char-poly}).
Therefore, the holonomy type for $R$ is the same as that of
$-\bar R$.   When Re$\,R\neq0$, we have two residues which give the same
holonomy. The two cases will be distinguished
by their vertical twist factors being $\infty$ or $-\infty$.

Also, all holonomy types in the table actually occur for some $R\in\co$.
We check that for any $\lambda_1,\lambda_2,\lambda_3\in\re$ with
$\sum \lambda_i=0$, then the polynomial
\begin{equation} \label{poly-eq}
\prod (\lambda-\lambda_i)= \chi(\lambda)=
\lambda^3-3(2^{-\frac23})|R|^{\frac23}\lambda
- {\rm Im }\,R,
\end{equation}
for some $R \in \co$.  We have to ensure that $|R|\ge |{\rm Im}\,R|$ with
equality only
in the case of multiple roots (by (\ref{disc-eq})).  Equation
(\ref{poly-eq}) is equivalent to
$$
{\rm Im}\,R=\lambda_1\lambda_2\lambda_3 \quad \mbox{and} \quad
-3(2^{-\frac23})|R|^{\frac23} =
\lambda_1\lambda_2+\lambda_1\lambda_3+\lambda_2\lambda_3
$$
for $\sum\lambda_i=0$.  Now write $\lambda_3=-\lambda_1-\lambda_2$.
By the homogeneity
properties of (\ref{poly-eq}) in $\lambda_i$ and $R$, we may assume
$\lambda_1=1$ also. We have
\begin{eqnarray*}
 |R|^2 -({\rm Im}\,R)^2 &=& \sfrac4{27}(\lambda_2^2+\lambda_2+1)^3-
(\lambda_2^2+\lambda_2)^2 \\
&=& \sfrac4{27}(\lambda_2+2)^2(\lambda_2-1)^2(\lambda_2+\sfrac12)^2
\end{eqnarray*}
This expression is always nonnegative and at each root,
$\lambda_i=\lambda_j$ for some $i\neq j$.

All together, we have

\begin{thm} \label{holonomy-thm}
On a Riemann surface $\Sigma= \bar\Sigma\setminus
\{p_i\}_{i=1\dots N}$, of negative Euler characteristic, and a
cubic form $U$ on $\Sigma$ which is allowed poles of order at most
3 at each puncture $p_i$, there is an $\rp^2$ structure.  The
$\rp^2$ holonomy at each end is determined by residue $R$ of $U$
at the corresponding puncture, i.e.\ by the $z^{-3}$ coefficient
in the Laurent series of $U$, by Table \ref{holonomy-table}.
Conversely, every hyperbolic, quasi-hyperbolic, and parabolic
holonomy type is determined by some $R\in\co$.
\end{thm}

\section{Detailed structure of the ends} \label{detail-sec}

\subsection{The triangle model} \label{triangle-model}
Recall the situation above. Model a neighborhood of a puncture of
$\Sigma$ by a punctured disc $\mathbb{D}_0$, and let the map
$\zeta(w)=e^{iw}$ be the covering map from the upper half-plane
$\mathbb H$ to $\mathbb{D}_0$.  Recall the asymptotics for $U$
(\ref{U-upstairs}). For our model metric $h$ (\ref{flat-metric}),
\begin{equation} \label{h-upstairs}
\zeta^*h=2^{\frac13}|R|^{\frac23}|dw|^2
\end{equation}
for $y \gg0$.  Then if we change coordinates
\begin{equation} \label{xi-eq}
\xi^3=2iR^{-1},\qquad \nu=\sigma+i\tau=\xi w,
\end{equation}
then $U=2\,d\nu^3$, $h=2\,|d\nu|^2$, and we can use this model
for any nonzero residue $R$.

So consider the complex plane $\mathbb C$ with coordinates
$\nu=\sigma+i\tau$, metric $h=2\,|d\nu|^2$ and cubic form
$U=2\,d\nu^3$. This configuration of $h$ and $U$ satisfies the
conditions above to form an affine sphere.  In fact, we can
explicitly solve the initial value problem.  From (\ref{fzz-eq})
for the frame $\{f,f_\nu,f_{\bar{\nu}}\}$
\begin{eqnarray*}
\frac{\partial}{\partial \sigma}\left(
\begin{array}{c}
f\\ f_\nu \\ f_{\bar \nu} \\
\end{array}
\right)
&=&
\left( \begin{array}{ccc}
0&1&1\\
1&0&1\\
1&1&0\\
\end{array} \right)
\left(
\begin{array}{c}
f \\ f_\nu \\ f_{\bar\nu} \\
\end{array}
\right),
 \\
\frac{\partial}{\partial \tau}\left(
\begin{array}{c}
f \\ f_\nu \\ f_{\bar\nu} \\
\end{array}
\right)
&=&
\left( \begin{array}{ccc}
0&i&-i\\
-i&0&i\\
i&-i&0\\
\end{array} \right)
\left(
\begin{array}{c}
f \\ f_\nu  \\ f_{\bar \nu} \\
\end{array}
\right).
\end{eqnarray*}
Since the two matrices above are simultaneously diagonalizable
(they must commute since the system is integrable), we can solve
this system explicitly to find that $(f,f_\nu,f_{\bar\nu})^\top$
is
\begin{equation} \label{explicit-sol}
\frac13
\left( \begin{array}{ccc}
1&1&1\\
1&\omega^2&\omega\\
1&\omega&\omega^2\\
\end{array} \right)
\left( \begin{array}{ccc}
e^{2\sigma}&0&0\\
0&e^{-\sigma+\sqrt{3}\tau}&0\\
0&0&e^{-\sigma-\sqrt{3}\tau}\\
\end{array} \right)
\left( \begin{array}{ccc}
1&1&1\\
1&\omega&\omega^2\\
1&\omega^2&\omega\\
\end{array} \right)
\left(
\begin{array}{c}
f(0)\\ f_\nu(0) \\ f_{\bar\nu}(0) \\
\end{array}
\right) ,
\end{equation}
where $\omega=e^{\frac{2\pi i}3}=-\frac12+i\frac{\sqrt{3}}2$.  The
imbedding
$f$ is real and therefore $f_\nu(0)=\overline{f_{\bar\nu}(0)}$.
Also we have the initial condition (\ref{init-val}) so that
$$ \det(f(0),f_\nu(0),f_{\bar\nu}(0))=i.$$
Choose initial conditions according to the eigenvectors of $\bf A,
\bf B$: Let
\begin{equation} \label{model-init-con}
\left(
\begin{array}{c}
f(0)\\ f_\nu(0) \\ f_{\bar\nu}(0) \\
\end{array}
\right)
= \frac1{\sqrt3}
\left( \begin{array}{ccc}
1&1&1\\
1&\omega^2&\omega\\
1&\omega&\omega^2\\
\end{array} \right)
. \end{equation}
The affine sphere for any other choice of initial data will simply
differ by a map in $\Sl{3}$.

Now we have an explicit formula for the imbedding $f$ of the affine
sphere into $\re^3$:
\begin{equation} \label{triangle-aff-sph}
f= \sfrac1{\sqrt3}
\left( e^{2\sigma}, e^{-\sigma+\sqrt{3}\tau},e^{-\sigma-\sqrt{3}\tau}
\right).
\end{equation}
This affine sphere is asymptotic to a cone over a triangle:  Let
$\mathcal T$ be the triangle with vertices $v_1=[1,0,0]$,
$v_2=[0,1,0]$, and $v_3=[0,0,1]$ and interior
$\left\{[1,x_2,x_3]:x_2,x_3>0\right\}$ in
homogeneous coordinates in $\rp^2$.  Then the affine sphere in
(\ref{triangle-aff-sph}) is asymptotic to the boundary of the cone
over $\mathcal T$ in $\re^3$, which is the first octant.

Consider any ray in the $(\sigma,\tau)$ plane approaching infinity.
First use (\ref{triangle-aff-sph}) to put the path on the affine
sphere in $\re^3$, and then project down to $\rp^2$.  Then the image
of the ray approaches the boundary of the triangle $\mathcal T$ in a
way that depends on the angle of the ray in the usual polar
coordinates $(\sigma,\tau)=(r\cos\theta,r\sin\theta)$. We record this in
Table \ref{limit-pts}, and note that in the cases where the limit point
is on a line segment, the exact
limit point is determined by the $\tau$-intercept of the ray.

\begin{table}
\caption{\label{limit-pts} Limit Points}
{
\renewcommand{\arraystretch}{1.3}
$
\begin{array}{|c|c|}
\hline
\mbox{range of }\theta& \mbox{limit point in }\partial{\mathcal T} \\
\hline \hline
(-\frac\pi3,\frac\pi3) & v_1 \\
\hline
\frac\pi3 & \mbox{on segment }\overline{v_1v_2}\\
\hline
(\frac\pi3,\pi) & v_2 \\
\hline
\pi & \mbox{on segment }\overline{v_2v_3} \\
\hline
(\pi,\frac{5\pi}3) & v_3 \\
\hline
\frac{5\pi}3 & \mbox{on segment }\overline{v_3v_1}\\
\hline
\end{array}
$
}
{}
\end{table}

Now pass back to the $w$ coordinate as in (\ref{xi-eq}).  The
topology of the end specifies two things. First, a clockwise
orientation around the puncture of the loop $|z|=\epsilon$ pulls
back to give the direction $\frac\partial{\partial x}$ with which
we have computed the holonomy. Now relate this holonomy direction
to $R$: Consider arg$\,R\in[-\frac\pi2,\frac{3\pi}2)$. For an
appropriate choice of cube root in (\ref{xi-eq}),
arg$\,\xi\in(-\frac\pi3,\frac\pi3]$. (Note the cube root needed to
find $\xi$ corresponds to the threefold symmetry of the triangle
$\mathcal{T}$.)  Then the holonomy direction in the $\nu$ plane is
equal to arg$\,\xi$. (We normalize the direction
$\frac\partial{\partial\sigma}$ in the $\nu$ plane to be $0$.)

Second, we have that any ray in the $w$ plane in any direction between
the $w$ plane of the form
\begin{equation} \label{dir-to-end}
\cos\iota\frac\partial{\partial x} +
\sin\iota\frac\partial{\partial y}
\end{equation}
for $\iota\in(0,\pi)$ will approach the end.  In the $\nu$ plane,
then, any direction between arg$\,\xi$ and arg$\,\xi+\pi$ will
approach the end. Below we will choose particular rays going to
the puncture to map out the affine sphere and determine the
vertical twist parameter for a puncture with residue $R$ if
Re$\,R\neq0$.

\begin{figure}
\begin{center}
\scalebox{.4}{\includegraphics{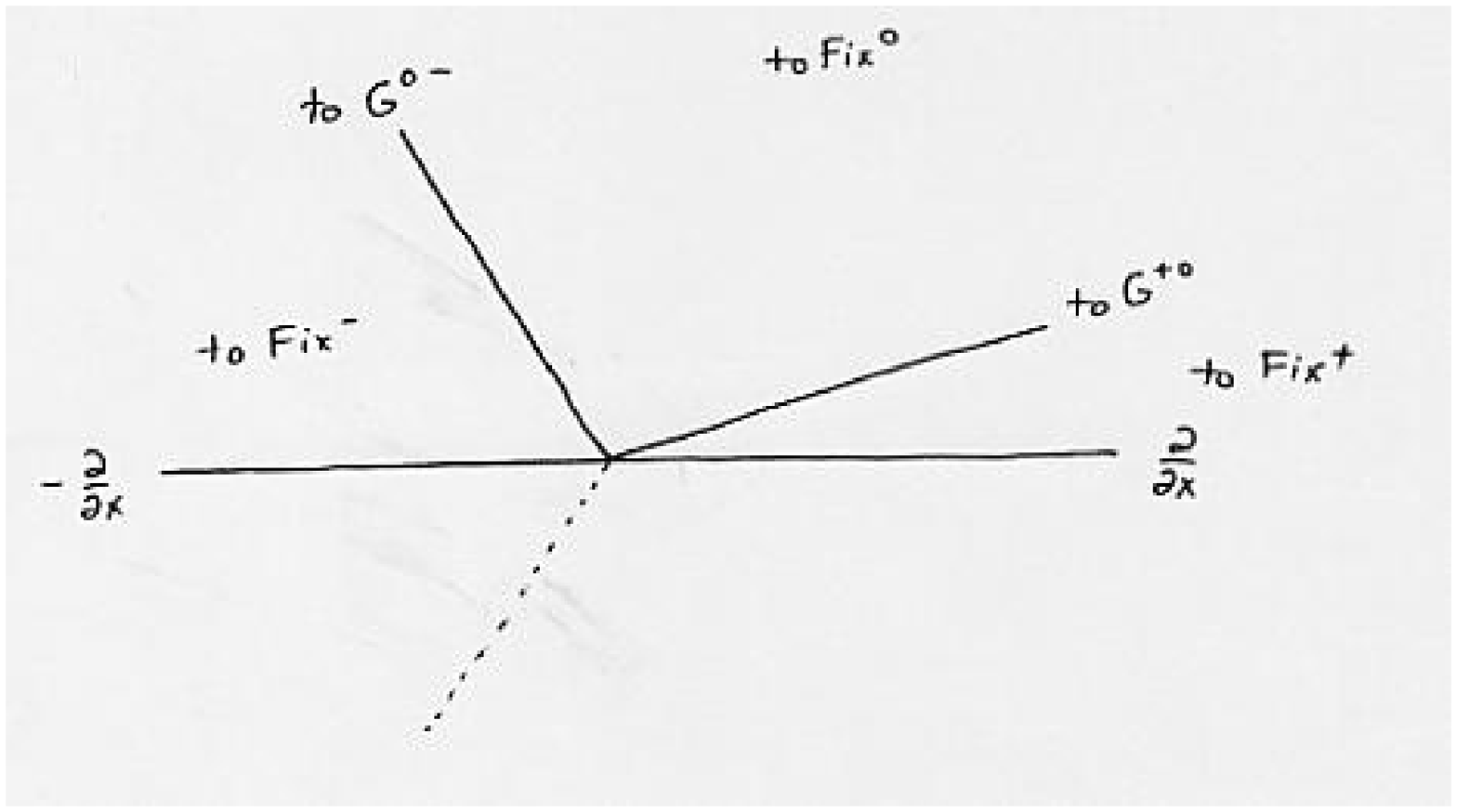}}
\end{center}
\caption{Rays in the $w$ Plane for Re$\,R>0$ \label{twist-pos}}
\end{figure}

Notice that three things can happen depending on the sign of
Re$\,R$.  If Re$\,R>0$, then arg$\,\xi\in(0,\frac\pi3)$.  See
Figure \ref{twist-pos}.  Then for $\theta\in({\rm arg}\,\xi,\frac\pi3)$, the
ray goes to the attracting fixed point of the holonomy, and if
$\theta\in(\pi,\pi+{\rm arg}\,\xi)$, the ray goes to the repelling
fixed point of the holonomy. (The limit points of the remaining
rays for $\theta=\frac\pi3,\pi$ should map out the geodesics
between the corresponding fixed points.) This leaves rays with
$\theta\in(\frac\pi3,\pi)$ to go to the saddle fixed point and
thus we should have the vertical twist parameter $\infty$.

\begin{figure}
\begin{center}
\scalebox{.4}{\includegraphics{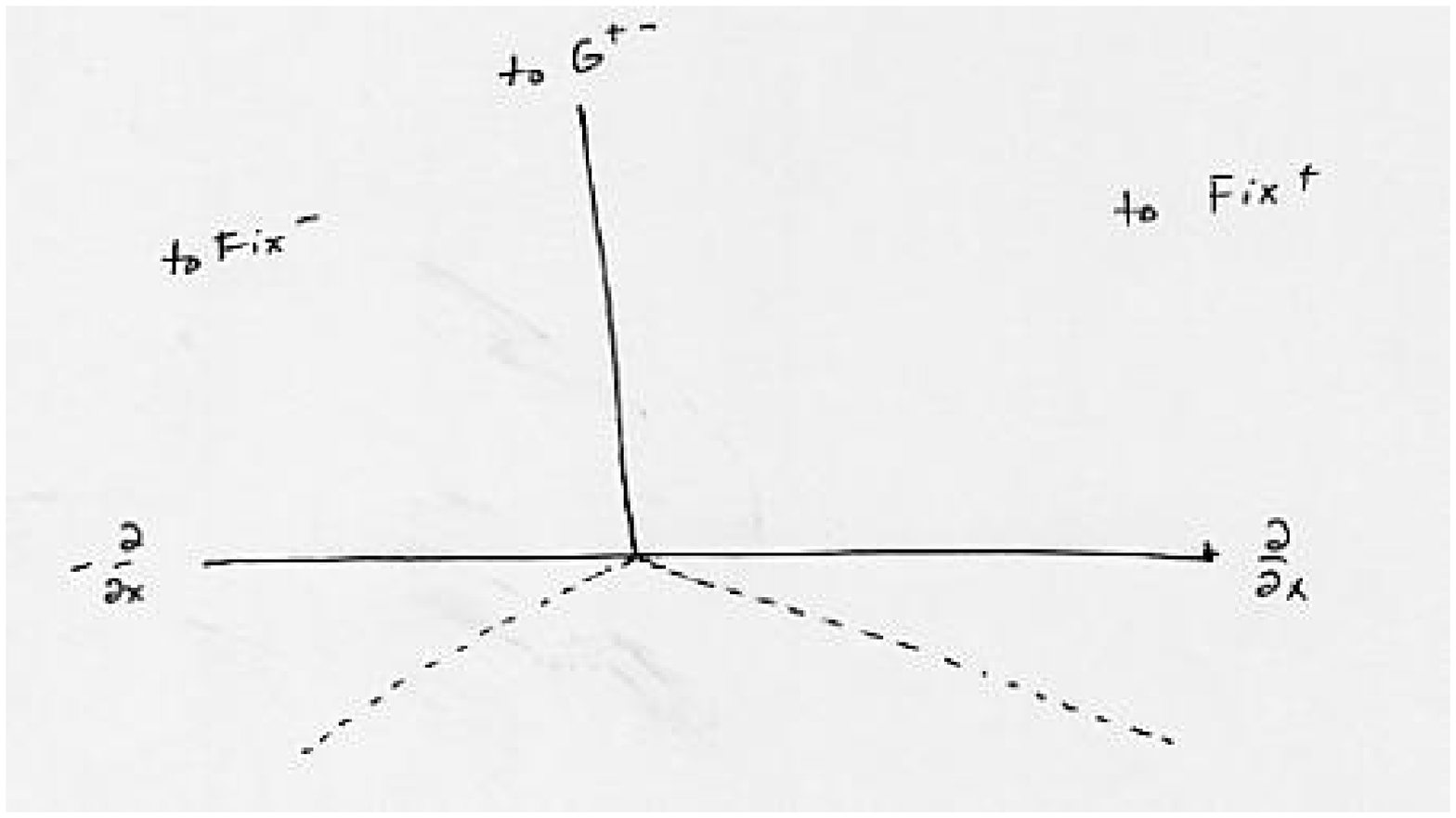}}
\end{center}
\caption{Rays in the $w$ Plane for Re$\,R<0$ \label{twist-neg}}
\end{figure}

On the other hand if Re$\,R<0$, arg$\,\xi\in(-\frac\pi3,0)$. See
Figure \ref{twist-neg}.  Then for $\theta\in({\rm arg}\,\xi,\frac\pi3)$, the
ray goes to the attracting fixed point of the holonomy, and if
$\theta\in(\frac\pi3,\pi+{\rm arg}\,\xi)$, the ray goes to the
repelling fixed point of the holonomy.  The limit points of the
ray with $\theta=\frac\pi3$, should map out the geodesic between
these two fixed points and the vertical twist parameter will be
$-\infty$.

If Re$\,R=0$, then arg$\,\xi=0$ or $\frac\pi3$, and the model
breaks down. It would predict holonomy $F$ as in Proposition
\ref{q-h-holo-prop} above, which we know is incorrect.

The next few subsections will prove the $\rp^2$ structure of an
end with Re$\,R\neq0$ follow the predictions we have just made.

\subsection{Perturbed linear systems} \label{perturb}
This model for the developing map is valid only near
a given puncture.  Recall the basic setup:  We lift
a neighborhood of a puncture $\{z\!: \epsilon>|z|>0\}$ on $\Sigma$ to the
region in the upper half plane $\{w=x+iy\!: y>-\log\epsilon\}$. In this
region we'll use our bounds on $u$ to approximate the initial value
problem
for the affine sphere (\ref{fzz-eq}) by the explicit models
computed above.  As discussed in the previous subsection, it may not be
useful only to consider the direction $\frac\partial{\partial y}$ going
to infinity, but also other directions of the form (\ref{dir-to-end}).
So for $\iota\in(0,\pi)$, introduce new coordinates
$\tilde x=x-y\,\cot\iota$,
$\tilde y= y\,\csc \iota$ so that
{
\renewcommand{\arraystretch}{2}
\begin{eqnarray} \label{xtilde-def}
\D \frac\partial{\partial \tilde x} &=& \D\frac\partial{\partial x} \\
\label{ytilde-def}
\D \frac\partial{\partial \tilde y} &=&
\D\cos\iota\,\frac\partial{\partial x}
                                 + \sin\iota\,\frac\partial{\partial y}
\end{eqnarray}
}

From (\ref{fzz-eq}) we have the equations for $\{f,f_w,f_{\bar w}\}$
\begin{eqnarray} \label{x-eq}
\frac{\partial}{\partial x}\left(
\begin{array}{c}
f\\ f_w \\ f_{\bar w} \\
\end{array}
\right)
&=&
\left( \begin{array}{ccc}
0&1&1\\
\frac12e^{\psi} & \psi_w& Ue^{-\psi} \\
\frac12e^{\psi} & \bar{U}e^{-\psi} & \psi_{\bar w}
\end{array} \right)
\left(
\begin{array}{c}
f \\ f_w \\ f_{\bar w} \\
\end{array}
\right),
 \\ \label{y-eq}
\frac{\partial}{\partial y}\left(
\begin{array}{c}
f \\ f_w \\ f_{\bar w} \\
\end{array}
\right)
&=&
\left( \begin{array}{ccc}
0&i&-i\\
-i\frac12e^\psi & i\psi_w &iUe^{-\psi}\\
i\frac12e^\psi & -i\bar{U}e^{-\psi} & -i\psi_{\bar w}\\
\end{array} \right)
\left(
\begin{array}{c}
f \\ f_w  \\ f_{\bar w} \\
\end{array}
\right),
\end{eqnarray}
and corresponding equations in the $\tilde x$ and $\tilde y$
coordinates. Recall the metric $e^\psi|dw|^2=e^uh$.  The bounds on
$u$ given in Proposition \ref{bound-prop} show that in the $w$
coordinates each of $u, \,u_w = O(e^{-2\alpha y})$ as $y\to\infty$
for some small positive $\alpha$.  Then along with the asymptotics
of $U$ in (\ref{U-upstairs}), the definition of $h$
(\ref{flat-metric}), and (\ref{xtilde-def}-\ref{ytilde-def}), we
have the asymptotic result
\begin{eqnarray}
\frac{\partial}{\partial \tilde x}\mathbf{X}&=&
\left[\mathbf{A}+\mathbf{O}(e^{-2\alpha y})\right]\mathbf{X}
\label{x-dir} \\
\frac{\partial}{\partial \tilde y}\mathbf{X}&=&
\left[\mathbf A \,\cos\iota + \mathbf B \,\sin \iota +
\mathbf O  (e^{-2\alpha y})\right]\mathbf{X}
\label{y-dir}
\end{eqnarray}
where $\mathbf{X}=(f,f_w,f_{\bar w})^{\top}$,
\begin{eqnarray*}
\mathbf{A} &=&
\left( \begin{array}{ccc}
0&1&1\\
2^{-\frac23}|R|^{\frac23} & 0 & -i 2^{-\frac13}R|R|^{-\frac23} \\
2^{-\frac23}|R|^{\frac23} & i2^{-\frac13}\bar{R}|R|^{-\frac23} & 0
\end{array} \right), \\
\mathbf{B} &=&
\left( \begin{array}{ccc}
0&i&-i\\
-i2^{-\frac23}|R|^{\frac23} & 0 & 2^{-\frac13}R|R|^{-\frac23} \\
i2^{-\frac23}|R|^{\frac23} & 2^{-\frac13}\bar{R}|R|^{-\frac23} & 0\\
\end{array} \right).
\end{eqnarray*}

In order to solve this integrable system, first solve the system
in the $\tilde y$ direction from some initial condition---this
will just be an ODE---and then solve in the $\tilde x$ direction
(or vice versa).  Consider the system
\begin{eqnarray}
\label{const-x-coeff}
\partial_{\tilde x}\mathbf{X} &=& \mathbf{AX}
\\
 \label{const-coeff}
\partial_{\tilde y}\mathbf{X} &=& (\cos\iota\,\bf A+\sin\iota\,\bf
B)\bf X.
\end{eqnarray}

Now change to the $\nu$ coordinate to relate these to the explicit
formulas in (\ref{explicit-sol}).  Then
\begin{eqnarray*}
 \frac\partial{\partial\tilde x} &=&
{\rm Re}\,\xi \,\frac\partial{\partial\sigma}
+ {\rm Im}\,\xi \,\frac\partial{\partial\tau}\\
\frac\partial{\partial\tilde y} &=&
{\rm Re}(\xi e^{i\iota}) \,\frac\partial{\partial\sigma}
+ {\rm Im}(\xi e^{i\iota}) \,\frac\partial{\partial\tau}
\end{eqnarray*}
Then the equations (\ref{const-x-coeff}--\ref{const-coeff}) become
\begin{eqnarray}
\label{const-x-tilde}
 \partial_{\tilde x}\mathbf{Z} &=& \tilde{\mathbf{A}} \mathbf{Z} \\
\label{const-y-tilde}
\partial_{\tilde y}\mathbf{Z} &=& \tilde{\mathbf{B}} \mathbf{Z},
\end{eqnarray}
where $\mathbf{Z}=(f,f_\nu,f_{\bar \nu})^\top$ and
\begin{eqnarray*}
\tilde{\mathbf{A}} &=&
\mathbf{P}\,\mathbf{D}\left(\rho_1,\rho_2,\rho_3
\right)
\,\mathbf{P}^{-1}, \\
\tilde{\mathbf{B}} &=&
\mathbf{P}\, \mathbf{D}\left( \mu_1,\mu_2,\mu_3\right)
\,\mathbf{P}^{-1},
\end{eqnarray*}
\begin{equation} \label{a-hat}
\begin{array}{l@{\qquad}l}
\mu_1 = 2 \,{\rm Re} (\xi e^{i\iota}), &
     \rho_1=2 \,{\rm Re}\,\xi, \\
\mu_2 = 2 \,{\rm Re} (\xi e^{i(\iota-\frac{2\pi}3)}), &
     \rho_2=2 \,{\rm Re} (\xi e^{-i\frac{2\pi}3}), \\
\mu_3 = 2 \,{\rm Re} (\xi e^{i(\iota+\frac{2\pi}3)}), &
     \rho_3=2 \,{\rm Re} (\xi e^{i\frac{2\pi}3}),
\end{array}
\end{equation}
\begin{equation}
\label{eigenvectors}
\mathbf{P}=
\frac1{\sqrt3}
\left( \begin{array}{ccc}
1 & 1 & 1 \\
1 & \omega & \omega^2   \\
1 & \omega^2 & \omega
\end{array} \right).
\end{equation}
Here $\mathbf D(\cdot,\cdot,\cdot)$ is a diagonal matrix.
Also (\ref{x-dir}--\ref{y-dir}) become
\begin{eqnarray}
\label{x-tilde-dir}
\partial_{\tilde x} \mathbf{Z} &=& \left[\tilde{\mathbf{A}} +
\mathbf{O}\left(e^{-2\alpha\tilde y\sin\iota}\right) \right]\mathbf{Z}\\
\label{y-tilde-dir}
\partial_{\tilde y}\mathbf{Z} &=& \left[\tilde{\mathbf{B}} +
\mathbf{O}\left(e^{-2\alpha\tilde y\sin\iota}\right) \right] \mathbf{Z}
\end{eqnarray}
Denote the perturbation terms $\mathbf{O}\left(e^{-2\alpha\tilde
y\sin\iota}\right)$ by $\mathbf{A}^{\rm pert}$ and
$\mathbf{B}^{\rm pert}$ respectively.

There is a theory, developed originally by Dunkel \cite{dunkel},
which addresses solutions to (\ref{y-tilde-dir}) as a perturbation
of (\ref{const-y-tilde}) as $\tilde y\to\infty$. Below in Appendix
\ref{app}, we follow Levinson \cite{levinson48} to find a version
which works with parameters. See also Hartman-Wintner
\cite{hartman-wintner55}. It is convenient to state it in terms of
an eigenbasis for $\tilde{\mathbf B} $.
\begin{lem} \label{pert-sol}
Let $Z_i$ be a vector solution of the linear constant-coefficient
equation (\ref{const-y-tilde}) which has the form $e^{\mu_i\tilde
y}\chi_i$ for a vector $\chi_i$.  Then (\ref{y-tilde-dir}) has a
solution $W_i$ such that
\begin{equation} \label{perturb-asymp}
\lim_{\tilde y\to\infty}\|W_i-Z_i\|e^{-\mu_i \tilde y}=0.
\end{equation}
Moreover, for any solution $W$ to (\ref{y-tilde-dir}), there is a
solution $Z$ to (\ref{const-y-tilde}) satisfying estimate
(\ref{perturb-asymp}) (with $\mu_i$ replaced by $\mu=\lim_{\tilde
y\to\infty}\frac1{\tilde y}\log\|W\|$).

If (\ref{y-tilde-dir}) depends continuously on a set of parameters
$\tau$---i.e.\ if the perturbation term $\,\mathbf{B}^{\rm
pert}(\tilde y, \tau)$ continuously varies in $C^0([T,\infty))$ as
$\tau$ varies---then $W_i(\tilde y,\tau)$ is continuous in $\tau$
and the limit (\ref{perturb-asymp}) is uniform in $\tau$.
\end{lem}
\begin{proof}
The exponential decay of the perturbation term in
(\ref{y-tilde-dir}) are more than sufficient to apply the results
in Levinson \cite{levinson48} to prove the first statement. The
second statement follows from the first by the fact that we can
choose $\{W_i\}$ to be a basis for the solution space of
(\ref{y-tilde-dir}).  We check in Appendix \ref{app} that the
estimates are uniform in a parameter $\tau$.
\end{proof}

Choose eigenvectors $\chi_i$ to be the column vectors of $\mathbf P$
(\ref{eigenvectors}).  Then by the lemma we have a matrix solution
$\mathbf Z$ to (\ref{y-tilde-dir}) so that as $\tilde y\to\infty$,
\begin{equation} \label{z-sol}
 \mathbf Z = \mathbf{P}
\left( \begin{array}{ccc}
e^{\mu_1\tilde y} + o(e^{\mu_1\tilde y}) &
o(e^{\mu_2\tilde y}) &
 o(e^{\mu_3\tilde y}) \\
 o(e^{\mu_1\tilde y}) &
 e^{\mu_2\tilde y} + o(e^{\mu_2\tilde y}) &
 o(e^{\mu_3\tilde y}) \\
 o(e^{\mu_1\tilde y}) &
 o(e^{\mu_2\tilde y}) &
 e^{\mu_3\tilde y}  + o(e^{\mu_3\tilde y})
\end{array} \right).
\end{equation}
This is enough to show that the limit of the ray in the direction
$\frac\partial{\partial \tilde y}$ approaches a point on the boundary
of the triangle $\mathcal T$, as in Table \ref{limit-pts}, with
direction $\theta=\iota+ {\rm arg}\,\xi.$ We also need to know,
however, how this solution (\ref{z-sol}) relates to the holonomy.
Equation (\ref{z-sol}) provides us with some information about a
frame at infinity and how to relate it to the frame
$\{f,f_\nu,f_{\bar\nu}\}$.

Pick a point $p_0$ so that $\tilde x =0$ and $\tilde y=\tilde y_0 \gg0$.
For any $\tilde y$, let consider the holonomy matrix $\mathbf{H}_{\tilde
y}$ for the frame $\{f,f_\nu,f_{\bar \nu}\}$ along the loop $(\tilde x,
\tilde y)$
for $\tilde x \in [0,2\pi]$.
Denote $\mathbf{H}_0=\mathbf{H}_{\tilde y_0}.$
By Section \ref{wang-corr} above, we are free to choose initial conditions
for $f$ in the initial value
problem (\ref{fzz-eq}), at least up to a multiplicative constant, which will
not essentially affect our arguments.
So choose initial conditions for $f,f_\nu,f_{\bar \nu}$ at $p_0$ so that for
$\mathbf Z$ from (\ref{z-sol})
\begin{equation} \label{init-z-cond}
\left(f(p_0),f_\nu(p_0), f_{\bar \nu}(p_0)\right)^\top =
\mathbf{Z}(p_0)=\mathbf{Z}_0
\end{equation}
Then in a path from $(0,\tilde y_0)$ to $(0,\tilde y)$, the holonomy
with respect to our frame is
\begin{equation} \label{y-tilde-hol}
\mathbf{F}_{\tilde y}= \mathbf{Z}_{\tilde y}
\mathbf{Z}_0^{-1}.
\end{equation}
Since the connection $D$ is flat,
\begin{equation} \label{flat-hol}
\mathbf{H}_{\tilde y} \mathbf{F}_{\tilde y} =
\mathbf{F}_{\tilde y} \mathbf{H}_0.
\end{equation}
We also know, as in (\ref{hol-limit}), $\lim_{\tilde y \to\infty}
\mathbf{H}_{\tilde y} = e^{2\pi\tilde{\mathbf{A}}}$, and so as
$\tilde y\to\infty$,
\begin{equation} \label{hol-limit2}
\mathbf{H}_{\tilde y} = \mathbf{P} \left[
\mathbf{D}(e^{2\pi\rho_1},e^{2\pi\rho_2},e^{2\pi\rho_3})
+ \mathbf{o}(1) \right] \mathbf{P}^{-1},
\end{equation}
where $\rho_i$ are the eigenvalues of $\tilde{\mathbf{A}}$ as in
(\ref{a-hat}).
Compute by (\ref{z-sol}), (\ref{y-tilde-hol}), (\ref{flat-hol}),
and (\ref{hol-limit2})
\begin{equation} \label{h-naught}
 \mathbf{H}_0= \mathbf{Z}_0
\left( \begin{array}{ccc}
e^{2\pi\rho_1} + o(1) &
o(e^{(-\mu_1+\mu_2)\tilde y}) &
 o(e^{(-\mu_1+\mu_3)\tilde y}) \\
 o(e^{(\mu_1-\mu_2)\tilde y}) &
 e^{2\pi\rho_2} + o(1) &
 o(e^{(-\mu_2+\mu_3)\tilde y}) \\
 o(e^{(\mu_1-\mu_3)\tilde y}) &
 o(e^{(\mu_2-\mu_3)\tilde y}) &
 e^{2\pi\rho_3}  + o(1)
\end{array} \right)
\mathbf{Z}_0^{-1}.
\end{equation}

In the applications below, we choose the parameter $\iota$
so that two of the eigenvalues $\mu_i$ of $\tilde{\mathbf{B}}$ are equal
to each other and greater than the remaining eigenvalue---see
(\ref{a-hat}). Assume without loss of generality
$\mu_1 = \mu_2 > \mu_3$. Then (\ref{hol-limit2}) shows
\begin{equation} \label{holonomy-in-frame}
\mathbf{H}_0= \mathbf{Z}_0
\left( \begin{array}{ccc}
e^{2\pi\rho_1}  &
0 &
 0 \\
 0 &
 e^{2\pi\rho_2}  &
 0 \\
 K &
 L &
 e^{2\pi\rho_3}
\end{array} \right)
\mathbf{Z}_0^{-1}.
\end{equation}
$K$ and $L$ denote real numbers over which
we have no a priori control, and this expression
is with respect to the frame $\mathbf{Z}_0=\{f(p_0),f_\nu(p_0),
f_{\bar \nu}(p_0)\}$.
Therefore, with respect to the standard frame in $\re^3$ then
the holonomy is given by
\begin{equation} \label{hol-eq-eig}
H=
\left( \begin{array}{ccc}
e^{2\pi\rho_1}  &
0 &
 0 \\
 0 &
 e^{2\pi\rho_2}  &
 0 \\
 K &
 L &
 e^{2\pi\rho_3}
\end{array} \right)
\end{equation}
This determines the eigenvalues of the holonomy (as in
Table \ref{holonomy-table}); also,
$(1,0,0)$ and $(0,1,0)$ are
eigenvectors corresponding to eigenvalues $e^{2\pi\rho_1}$ and
$e^{2\pi\rho_2}$ respectively.
(Note that $H$ acts on the right on row vectors in $\re^3$.)
Therefore,
projecting down to $\rp^2$, $[1,0,0]$ and $[0,1,0]$
are fixed points of the
holonomy action.  Each of these two fixed points
is attracting, saddle, or repelling
according to whether the appropriate $\rho_j$ ($j=1$ or 2)
is numerically the largest, the middle,
or the smallest of the $\{\rho_i\}_{i=1}^3$.

Recall that $f(0,\tilde y)$ is
the top row of the matrix $\mathbf Z$. Therefore, (\ref{z-sol}) shows
that, upon projecting from $\re^3$ to $\rp^2$,
$\lim_{\tilde y\to\infty}[f(0,\tilde y)]=[1,1,0]$.  Moreover, we have by
(\ref{hol-eq-eig})
$$ \lim_{\tilde y\to\infty}[f(2\pi n,\tilde y)]=
[e^{2\pi n\rho_1},e^{2\pi n\rho_2},0]=\ell_n,$$ for $n\in\mathbb
Z$. So all these limit points $\ell_n$ are on a geodesic segment
between $[1,0,0]$ and $[0,1,0]$.  Since they are limit points of
rays which go to infinity in the universal cover of $\Sigma$, and
since the $\rp^2$ structure is convex, these limit points $\ell_n$
in $\rp^2$ must all be on the boundary $\partial \Omega$.
Therefore, by convexity, the entire geodesic segment
$$\{[1,t,0]:t\in[0,\infty]\}\subset \partial\Omega.$$
We record this discussion in

\begin{prop} \label{fixed-point-type}
Let $v_i$ be the standard $i^{\rm th}$ basis vector in $\re^3$, and
$[v_i]$ the projection to $\rp^2$.
Choose the parameter $\iota$ so that
$\mu_j=\mu_k$ are the largest two eigenvalues of $\tilde{\mathbf B}$.
The points $[v_j]$ and $[v_k]$ are fixed points
of the holonomy.  Each is attracting, repelling or saddle-type according to
whether the corresponding eigenvalue of the holonomy
is numerically the largest, the smallest, or the
middle among $\{e^{2\pi\rho_i}\}_{i=1}^3$.  The line segment
$$ \{ [v_j+tv_k]: t\in(0,\infty) \}$$
is in the boundary of the image of the developing map.
\end{prop}

\subsection{Hyperbolic ends: the case ${\rm Re}\,R>0$}
\label{full-triangle-sec} In this case the vertical twist
parameter is $\infty$.

\begin{prop} \label{plus-prop}
Consider the end of the $\rp^2$ surface corresponding to
$(\Sigma,U)$ with residue $R$ at the end.  If Re$\,R>0$, then the
vertical twist parameter of the end is $\infty$.
\end{prop}

\begin{proof}
If Re$\,R>0$, then choose $\xi=(2iR^{-1})^\frac13$ so that
arg$\,\xi\in(0,\frac\pi3)$.  Let $\iota=\frac\pi3-{\rm arg}\,\xi$
and $\hat \iota= \pi-{\rm arg}\,\xi$ so that $\theta=\iota+{\rm
arg}\, \xi=\frac\pi3,$ and $\hat \theta=\pi$.

First consider the case $\theta=\frac\pi3$.
Then (\ref{a-hat}) shows that $\mu_1=\mu_2>\mu_3$ and
$\rho_1>\rho_2>\rho_3$.
For a value of $y\gg0$, choose initial condition (\ref{init-z-cond})
for the equation (\ref{fzz-eq}).
Proposition \ref{fixed-point-type}
shows that $[v_1]={\rm Fix}^+$,
$[v_2]={\rm Fix}^0$ and the geodesic segment $G^{+0}$ is contained
in the boundary of the image of the developing map.
Moreover, the holonomy matrix with respect
to the standard frame in $\re^3$ is given by (\ref{hol-eq-eig}).
Note the there is as yet no a priori control over
$K$ and $L$.  By measuring the holonomy in the $\hat\iota$
direction as well, we shall see that $K$ and $L$ do vary continuously
in families however.

Now for $\hat\iota$, $\hat\mu_2=\hat\mu_3>\hat\mu_1$ and we still have
$\rho_1>\rho_2>\rho_3$.  Therefore, again choose
choose initial condition (\ref{init-z-cond})
for the equation (\ref{fzz-eq}).  Note that this amounts to choosing
a new frame on $\re^3$,
$\{\hat v_i\}_{i=1}^3.$  With respect to this frame,
$[\hat v_2]={\rm Fix}^0$ and $[\hat v_3]={\rm Fix}^-$.  Proposition
\ref{fixed-point-type} shows that the geodesic segment $G^{0-}$
is in the boundary of the developing map.

By convexity, then the principal geodesic segment $G^{+-}$ must be contained
in the closure $\bar \Omega$ of $\Omega$ the image of the developing map.
We claim the open segment $G^{+-}\subset\Omega$.  If on the contrary
$G^{+-}\subset\partial\Omega$, then $\Omega$ is a triangle and as in
Subsection \ref{triangle-model},
the affine metric $h$ on $\Omega$ is complete and flat.
Therefore, $(\Omega,h)$ is conformally equivalent to $\co$.  This
contradicts the fact that $\Sigma$ admits a complete hyperbolic metric.
Now if there were a single point of the open segment $G^{+-}$ in
$\partial\Omega$, then since the endpoints ${\rm Fix}^+$, ${\rm Fix}^-
\in \partial \Omega$, convexity forces all of $G^{+-}\subset\partial
\Omega$, and we reach a contradiction again to prove the claim.

The discussion in Section \ref{goldman-coor} above then proves the
proposition.
\end{proof}

It will be useful below to show that the holonomy matrix (\ref{hol-eq-eig})
varies continuously in families.  We will determine the constants
$K$ and $L$ in terms of the change of frame between the $v_i$
and $\hat v_i$ above.

Diagonalize (\ref{hol-eq-eig})
$$ H=
\left( \begin{array}{ccc}
1 &
0 &
 0 \\
 0 &
 1  &
 0 \\
 K' &
 L' &
 1
\end{array} \right)
\left( \begin{array}{ccc}
e^{2\pi\rho_1}  &
0 &
 0 \\
 0 &
 e^{2\pi\rho_2}  &
 0 \\
 0 &
 0 &
 e^{2\pi\rho_3}
\end{array} \right)
\left( \begin{array}{ccc}
1  &
0 &
 0 \\
 0 &
 1  &
 0 \\
 -K' &
 -L' &
 1
\end{array} \right),
$$
where $K'=K/(e^{2\pi\rho_1}-e^{2\pi\rho_3})$ and
$L'=L/(e^{2\pi\rho_2}-e^{2\pi\rho_3})$.
Then in terms of the $\hat v_i$ frame,
\begin{eqnarray*}
\hat H &=&
 \left( \begin{array}{ccc}
e^{2\pi\rho_1}  &
\hat K &
 \hat L \\
 0 &
 e^{2\pi\rho_2}  &
 0 \\
 0 &
 0 &
 e^{2\pi\rho_3}
\end{array} \right) \\
&=&
\left( \begin{array}{ccc}
1 &
\hat K' &
 \hat L' \\
 0 &
 1  &
 0 \\
 0 &
 0 &
 1
\end{array} \right)
\left( \begin{array}{ccc}
e^{2\pi\rho_1}  &
0 &
 0 \\
 0 &
 e^{2\pi\rho_2}  &
 0 \\
 0 &
 0 &
 e^{2\pi\rho_3}
\end{array} \right)
\left( \begin{array}{ccc}
1  &
-\hat K' &
 -\hat L'\\
 0 &
 1  &
 0 \\
 0 &
 0 &
 1
\end{array} \right),
\end{eqnarray*}
where $\hat K'=\hat K/(e^{2\pi\rho_2}-e^{2\pi\rho_1})$ and
$\hat L'=\hat L/(e^{2\pi\rho_3}-e^{2\pi\rho_1})$.

Denote by $Q$ the change of frame in $\re^3$
between the $v_i$ and the $\hat v_i$.  Then $H=Q\hat H Q^{-1}$.
Since the eigenvalues $\{e^{2\pi\rho_i}\}$ are distinct, there must
exist real constants $\omega_i\neq0$ so that
$$
\left( \begin{array}{ccc}
1 &
0 &
 0 \\
 0 &
 1  &
 0 \\
 K' &
 L' &
 1
\end{array} \right)
\left( \begin{array}{ccc}
\omega_1 &
0 &
 0 \\
 0 &
 \omega_2  &
 0 \\
 0 &
 0 &
 \omega_3
\end{array} \right)
= Q
\left( \begin{array}{ccc}
1 &
\hat K' &
 \hat L' \\
 0 &
 1  &
 0 \\
 0 &
 0 &
 1
\end{array} \right)
$$
Let $Q=(Q^i_j)$.  The previous equation forces $Q^2_1=Q^2_3=0$, and the
free parameters are determined by $Q$:
$$
\begin{array}{c}
\D \hat K'=-\frac{Q^1_2}{Q^1_1}, \quad
\hat L'=-\frac{Q^1_3}{Q^1_1}, \quad
K'=\frac{Q^3_1}{Q^1_1}, \quad
L'=\frac{Q^3_2 Q^1_1-Q^3_1 Q^1_2}{Q^1_1 Q^2_2}, \\
\D \omega_1=Q^1_1, \quad
\omega_2=Q^2_2, \quad
\omega_3=\frac{Q^3_3 Q^1_1 - Q^3_1 Q^1_3}{Q^1_1}.
\end{array}
$$
Note the last row implies $Q^1_1,Q^2_2\neq0$. $Q$ is determined by
the initial conditions $\mathbf{Z}_0$ and $\hat{\mathbf{Z}}_0$
from (\ref{init-z-cond}) to the initial value problem
(\ref{fzz-eq}). Explicitly, Proposition \ref{two-sphere} shows
that $Q=\hat{\mathbf{Z}}_0^{-1} \mathbf{Z}_0$. $\mathbf{Z}_0$ and
$\hat{\mathbf{Z}}_0$ come from Lemma \ref{pert-sol}.  Therefore,
Lemma \ref{pert-sol} and Proposition \ref{app-prop} imply the
parameters $K,L,\hat K,\hat L$ vary continously in families.

\begin{prop}  \label{holo-matrix-param}
The holonomy matrix (\ref{hol-eq-eig}) constructed above varies
continuously for families of equations as long as the perturbation
term $\mathbf{B}^{\rm pert}(y,\tau)$ satisfies the hypotheses of
Lemma \ref{pert-sol}.
\end{prop}

\begin{rem}
Here is a more geometric interpretation of the preceeding Proposition and its
proof. The proposition allows us to control the $\rp^2$
coordinates of the developing map of a degenerating family.  In
other words, it allows us to control the gauge.  Developing along
rays for the parameter $\iota$ allow us to control the line segment
$G^{+0}$, and similarly for the parameter $\hat \iota$, we control
the line segment $G^{0-}$, up to a change of gauge $Q$ which varies
continuously.  Any automorphism of $\rp^2$ which fixes $G^{+0}\cup G^{0-}$
must be trivial, and so the proposition follows.
\end{rem}

\subsection{Hyperbolic ends: the case ${\rm Re}\,R<0$}
In this case, the vertical twist parameter is $-\infty$.

\begin{prop} \label{minus-prop}
On an end of $(\Sigma,U)$ with residue $R$ so that Re$\,R<0$,
the principal geodesic line segment of the holonomy action
around this end is in the boundary of the image of the developing map.
The vertical twist parameter for this end is $-\infty$.
\end{prop}

\begin{proof}
If Re$\,R<0$, then choose $\xi=(2iR^{-1})^\frac13$ so that
arg$\,\xi\in(-\frac\pi3,0)$.  Let $\iota=\frac\pi3-{\rm arg}\,\xi$
so that $\theta=\frac\pi3.$  Then (\ref{a-hat}) shows that
$\mu_1=\mu_2>\mu_3$ and
$\rho_1>\rho_3>\rho_2$.
Proposition \ref{fixed-point-type}
then implies that $[1,0,0]={\rm Fix}^+$ and $[0,1,0]={\rm Fix}^-$, and
that the geodesic segment between them, $G^{+-}$, is contained in the
boundary of the image of the developing map.
The discussion in Section \ref{goldman-coor} above then implies the vertical
twist parameter is $-\infty$.
\end{proof}

\subsection{Quasi-hyperbolic ends} \label{quasi-h-subsec}
For Re$\,R=0$, $R\neq0$, Table \ref{holonomy-table} shows the holonomy type
of the end is quasi-hyperbolic.  This completely determines the structure
of the end.

\begin{prop}
Let $S$ be a properly convex $\rp^2$ surface with an end with
quasi-hyperbolic holonomy.  If $\Omega \subset\rp^2$ is the universal
cover, then the end of $S$ has boundary given by the push-down of a geodesic
segment in $\partial \Omega$
whose endpoints are the two fixed points of the holonomy action.
\end{prop}

\begin{figure}
\begin{center}
\scalebox{.5}{\includegraphics{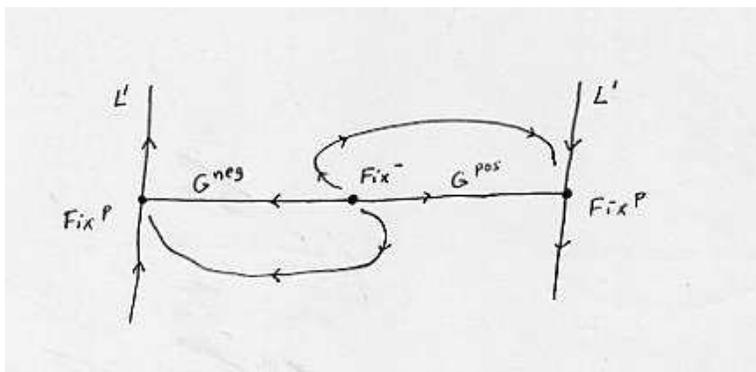}}
\end{center}
\caption{The Dynamics of Quasi-Hyperbolic Holonomy \label{quasihyperbolic-pic}}
\end{figure}

\begin{proof}
Lift the holonomy action to $\Sl{3}$ and choose coordinates
$(x_1,x_2,x_3)$ in $\re^3$ so that the holonomy matrix is of the
form
$$ \gamma=\left(
\begin{array}{ccc}
\alpha_1 & 1 & 0 \\
0 & \alpha_1 &0 \\
0 & 0 & \alpha_3
\end{array}
\right). $$ Here $\alpha_i>0$, $\alpha_1^2\alpha^3=1$, and
$\alpha_1\neq\alpha_3$. Consider the case $\alpha_1 > \alpha_3$.
The two fixed points of the holonomy are a repelling fixed point
Fix$^-=[0,0,1]$, and Fix$^P=[1,0,0]$. There are also two geodesic
lines preserved by the holonomy given by $L=\{x_2=0\}$, which
connects the two fixed points and on which the action is
hyperbolic; and $L'=\{x_3=0\}$, on which the action is parabolic.

Pick a point $p=[x_1,x_2,1]\in\Omega\setminus (L\cup L')$. Then
$p_n=\gamma^np$ satisfies $\lim_{n\to\infty}p_n= {\rm Fix}^P$,
$\lim_{n\to-\infty}p_n={\rm Fix}^-$.  The convexity of $\Omega$
implies that one of the two geodesic segments between Fix$^-$ and
Fix$^P$ must be contained in $\bar \Omega$.  If $x_2>0$, then it
is the segment $G^{\rm pos} =\{[t,0,1]:t>0\}$, and if $x_2<0$,
then it is the segment $G^{\rm neg}=\{[t,0,1]:t<0\}$.  Now we
claim that this segment $G^{\rm pos}$ or $G^{\rm neg}$ must be in
the boundary $\partial \Omega$. Without loss of generality, assume
$x_2>0$.  Then if a neighborhood of any point in $G^{\rm pos}$ is
contained in $\Omega$, then $\Omega$ contains a point with
$x_2<0$, and therefore $\bar \Omega$ contains $G^{\rm neg}$ as
well, and so $\partial \Omega$ contains the whole line $L$.  This
contradicts the fact that $\Omega$ is strictly convex.  See Figure
\ref{quasihyperbolic-pic}. The case
$\alpha_1<\alpha_3$ is similar.
\end{proof}

\subsection{Parabolic ends} \label{parabolic-subsec}
If the residue $R=0$, Table \ref{holonomy-table} shows the holonomy type
of the end is parabolic.  As in the quasi-hyperbolic case above, the
holonomy's being parabolic
completely determines the structure of the end.

\begin{prop}
Let $S$ be a properly convex $\rp^2$ surface with an end with
parabolic holonomy.  Assume the fundamental group $\pi_1(S)\neq\mathbb Z$.
If $\Omega \subset\rp^2$ is the universal
cover of $S$, then the end of $S$ has boundary (as a set)
given by the push-down of the single fixed point of the holonomy action,
which is in $\partial \Omega$.
\end{prop}

\begin{proof}
Lift the holonomy action to $\Sl{3}$.  Choose coordinates in
$(x^1,x^2,x^3)$ of $\re^3$ so that the holonomy matrix is
$$\gamma=\left( \begin{array}{ccc}
1&1&0 \\ 0&1&1 \\ 0&0&1
\end{array} \right).$$
The only fixed point is ${\rm Fix}=[1,0,0]$. Also, there is a single
line preserved
by the holonomy action given by the line $L=\{x^3=0\}$. The holonomy is
parabolic along $L$.  For any point $p\in \Omega$ and $p_n=\gamma^np$, then
we have $p_n\to {\rm Fix}$ for $n\to\pm\infty$.  Therefore,
$\rm{Fix}\in \bar\Omega$.  Since the holonomy acts on $\Omega$
without fixed points, $\rm{Fix}\in \partial \Omega$.

For a domain $\Omega$ on which $\gamma$ acts, there are two ends of
the cylindrical quotient $\Omega/\langle\gamma\rangle$. See Figure
\ref{parabolic-pic}.
Call the end which develop to Fix the small end, and the other the large
end.  We want to rule out the large end.  Notice that the holonomy action
$\gamma^n$ takes the part of $\partial \Omega$ associated with the large
end to all of $\partial \Omega \setminus \{ {\rm Fix} \}$.

\begin{figure}
\begin{center}
\scalebox{.5}{\includegraphics{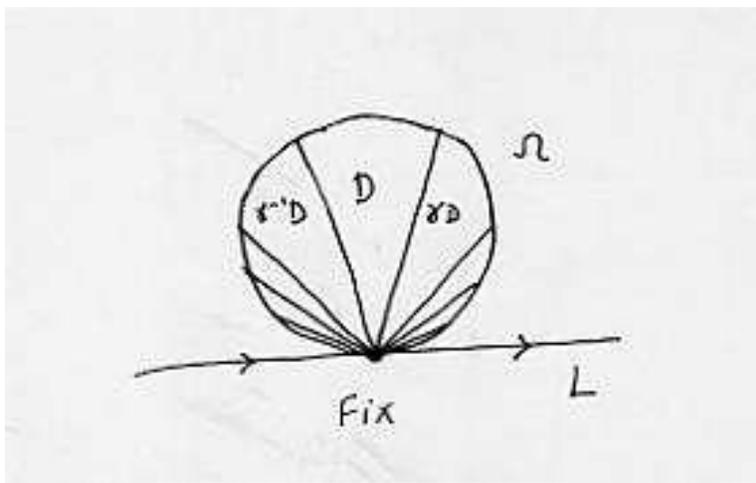}}
\end{center}
\caption{Parabolic Holonomy \label{parabolic-pic}}
\end{figure}

Assume that our end develops to the large end of $\Omega/\langle\gamma\rangle.$
As in Subsection \ref{setup} above, fix a basepoint $p\in S$ near the
end in question, and consider a path $P$
going from $p$ to the end which doesn't leave
a fixed cylindrical neighborhood $N$ of the end.  Fix a fundamental domain
$D\in\tilde S$ which contains a lift $\tilde P$ of the path $P$.  We may
assume that dev$(\tilde P)$ has a limit point $r\in \partial \Omega$. Since
we have dev$(\tilde P)$ is contained in the large end, we may assume
$r\neq{\rm Fix}$.  By the action of $\gamma^n$ on the large end, a simple
continuity argument, along with Figure \ref{parabolic-pic},
shows that the closure in $\rp^2$ of
$\bigcup_n{\rm dev}(\gamma^n D)$ contains all of $\partial \Omega$, and also
includes a neighborhood in $\bar \Omega$ of $\partial\Omega\setminus
\{ {\rm Fix} \}$.

Now consider an element
$\gamma'\in\pi_1(S)\setminus\langle\gamma\rangle$, and the
corresponding holonomy $H = {\rm hol}(\gamma')\in\pgl{3}$. Then we
may assume, by perturbing the path $P$ if necessary, that
$H(r)\neq {\rm Fix}$.  But of course $H (r)\in\partial \Omega$.
Also, $H (r)$ is the limit point of the path $\gamma'\tilde P$,
where we consider $\gamma'$ as the path starting at $\tilde p \in
\tilde S$ and covering an appropriate loop.  But this then
contradicts the fact that $\overline{\bigcup_n{\rm dev}(\gamma^n
D)}$ contains a neighborhood in $\bar \Omega$ of $\partial \Omega
\setminus \{ {\rm Fix} \}$.
\end{proof}

\section{The boundary of the moduli space} \label{deg-fam-sec}

\subsection{Regular 3-differentials}
The material in this subsection is well known.  I would like to thank
Michael Thaddeus, Ravi Vakil and Richard Wentworth for explaining some
of it to me.  A basic reference for the algebraic theory is the book
of Harris and Morrison \cite{harris-morrison}.  A good summary
of the analytic techniques used here is contained in Wolpert
\cite{wolpert02}. We only give a sketch of the arguments.

Consider the Deligne-Mumford compactification $\overline{\mathcal
M}_g$ of the moduli space of Riemann surfaces of genus $g\ge2$ and
also $\overline{\mathcal M}_{g,1}$ the compactification of the of
the moduli space of genus-$g$ Riemann surfaces with one marked
point $p$.  Let $\mathcal{F}\!: \overline{\mathcal M}_{g,1}\to
\overline{\mathcal M}_g$ be the forgetful map.  In this context,
$\overline{\mathcal M}_{g,1}$ is the \emph{universal curve} over
$\overline{\mathcal M}_g$. These moduli spaces are only
V-manifolds (smooth Deligne-Mumford stacks). We will describe
complex coordinates on the sense of V-manifolds: in general, for
each point $x$ in the moduli space there is a chart given by the
quotient of an open set in $\co^N$ by a finite group $G$ of
biholomorphisms.  (For details of V-manifolds, see Baily
\cite{baily57}.) Charts in Teichm\"uller space provide these
V-manifold charts near $x\in\mathcal M_g$. (We describe below
coordinate charts near $x\in\overline{\mathcal M}_g \setminus
\mathcal M_g$.) The group $G$ is given by the group of
automorphisms of the curve $\mathcal{F}^{-1}(x)$ over the point
$x$.  (In the case $g=2$ with no marked point, the group $G$ is
instead the quotient of the group of automorphisms of the curve by
the automorphism group of the generic curve, which is generated by
the hyperelliptic involution.)

At a point $x \in \mathcal M_g$, consider the curve $C_0$ over
$x$. In Teichm\"uller space $\mathcal T_g$ consisting of $3g-3$
Beltrami differentials $\nu_i$.  For $s=(s^i)\in\co^{3g-3}$ small,
then $\|s^i\nu_i\|_{L^\infty}<1$ and thus there is a
quasiconformal map homeomorphism $C_0$ to $C_s$ with Beltrami
differential $s^i\nu_i$.  Then $s$ form local $C^\infty$
V-manifold coordinates around $x\in\mathcal M_g$.  Note $s$ is
\emph{not} a holomorphic coordinate system.

Each point $x\in\overline{\mathcal M}_g\setminus \mathcal M_g$,
represents a Riemann surface with nodes.  In other words, at each
point $p$ in our curve $C$ over $x$ has a neighborhood of the form
either $\{z:|z|<K\}$ or $\{(z,z'):zz'=0,|z|<K,|z'|<K\}$. Let $n$
be the number of nodes.  Let $C^{\rm reg}$ be the smooth part of
$C$, which is formed by removing the $n$ nodes. $C^{\rm reg}$ is a
possibly disconnected noncompact Riemann surface. $C^{\rm reg}$
may be smoothly compactified to $\tilde C$ by adding $2n$ points
$\{p_i,q_i\}$. $\tilde C$ is the \emph{normalization} of $C$. A
natural analytic map from $\tilde C$ to $C$ identifies each pair
$(p_i,q_i)$ to a single point to form each node.  $C$ must be
stable: i.e.\ $C$ has only finitely many automorphisms;
equivalently $C^{\rm reg}$ admits a conformal complete hyperbolic
metric.

A natural deformation of $C$ consists of \emph{plumbing} the
nodes, in which we replace each node $\{zz'=0\}$ by the smooth
neck $\{zz'=t\}$.  There is a nice overview of the plumbing
construction in Wolpert \cite{wolpert02}.  A neighborhood of the
$i^{\rm th}$ node is for local coordinates $z_i,z_i'$
\begin{equation}\label{node-neigh}
N_i=\{(z_i,z_i'):z_iz_i'=0,|z_i|<K,|z_i'|<K\}.
\end{equation}
For a complex parameter $t_i:|t_i|<K^2$, we replace $N_i$ by the
smooth cylinder
$N_i^{t_i}=\{(z_i,z_i'):z_iz_i'=t_i,|z_i|<K,|z_i'|<K\}$. It is
clear that for $t$=$(t_i)$
$$C_t=\left[C \setminus \left(\cup_i N_i\right)\right] \bigsqcup
\left(\cup_i N_i^{t_i}\right)$$ patches together
complex-analytically to make $\Sigma_t$ smooth on a neighborhood
of each $\overline{N_i^{t_i}}$.  These $t_i$ form $n$ complex
V-manifold coordinates over $\overline{\mathcal M}_g$.

In addition, Wolpert \cite[Lemma 1.1]{wolpert92b} has shown there
is a real-analytic family of Beltrami differentials $\nu(s)$ on
$C^{\rm reg}$  parametrized by $s$ in a neighborhood of the origin
in $\co^{3g-3-n}$ so that the induced quasiconformal maps
$\zeta^{\nu(s)}\!:C^{\rm reg}\to C^{\rm reg}_s$ preserve the cusp
coordinate up to multiplying by a rotation $e^{i\theta_s}$.  Form
$C_s$ by completing $C^{\rm reg}_s$ by reattaching the
corresponding nodes. For the $n$ nodes and $t\in\co^n$ small, we
perform the plumbing construction as above with respect to the
cusp coordinates on each $C_s$ to form a family $C_{s,t}$ of
curves. Then $(s,t)$ form a real-analytic V-manifold coordinate
chart of $\overline{\mathcal M}_g$ near a nodal curve.  For each
fixed $s$, the $t$ coordinates are complex-analytic.

There are similar coordinates on the universal curve
$\overline{\mathcal M}_{g,1}$.  The idea is to treat the marked
point $p$ as a puncture.  As above, by Lemma 1.1 in
\cite{wolpert92b}, in $\mathcal M_{g,1}$, there is a real-analytic
family of Beltrami differentials $\nu(s)$ for $s$ in a
neighborhood of 0 in $\co^{3g-3}$ so that the quasiconformal maps
$\zeta^{\nu(s)}$ preserve a canonical complex coordinate
neighborhood $U$ of $p\in C_0$.  Then $p$ may move
complex-analytically in $U$ so that $(s,p)$ form a real-analytic
V-manifold coordinate chart of $\mathcal M_{g,1}$, and for each
fixed $s$, the $p$ coordinate is complex-analytic.

It is straightforward to combine the construction in the last two
paragraphs in the case of a nodal curve $C$ with $n$ nodes over a
point in $\overline{\mathcal M}_g \setminus \mathcal{M}_g$ and a
point $p\in C^{\rm reg}$.   Then we have real-analytic coordinates
$(s,p,t)$ with $s\in U\subset \co^{3g-3-n}$, $t\in U' \subset
\co^n$, and $p\in U'' \subset C^{\rm reg}_{s,t}$ so that for $s$
fixed, the coordinates $(p,t)$ are complex-analytic.

\begin{figure}
\begin{center}
\scalebox{.5}{\includegraphics{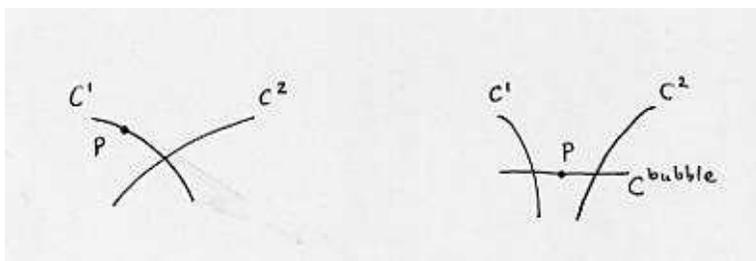}}
\end{center}
\caption{$\overline{\mathcal M}_{g,1}$ near a Node \label{C-bubble-pic}}
\end{figure}

In the remaining case of a nodal curve $C$ with $p$ at a node is
more subtle.  First of all, having $p$ equal to a node is
technically not allowed in the Deligne-Mumford compactification.
Let the node in $C$ be represented by $zz'=0$.  Then if $p$ is the
point of the node $\{z=z'=0\}$, this configuration is not stable.
Instead introduce a sphere $\mathbb{CP}^1$ attached to $\tilde C$
by one point $r\in\mathbb{CP}^1$ to $\{z=0\}$ and by another point
$r'\in\mathbb{CP}^1\setminus\{r\}$ to $\{z'=0\}$.  This amounts to
having a sphere ``bubbling" to separate the existing node into two
pieces and having the sphere attached by nodes at $r,r'$ to each
piece. See Figure \ref{C-bubble-pic} (in the case the node
separates the curve into two parts $C^1$ and $C^2$).
Then $p$ is allowed to be any point in $\mathbb{CP}^1
\setminus \{r,r'\}$. However, since the sphere with three marked
points $p,r,r'$ has no automorphisms, we may collapse the
$\mathbb{CP}^1$ and identify this configuration canonically with
the configuration of the point $p$ equal to the node $\{z=z'=0\}$.
Then $p$ is no longer a smooth complex coordinate, since it must
vary in a singular curve at the singularity. Instead consider the
\emph{plumbing variety} $\{zz'=t\}\subset \co^3$ for $z,z',t$ near
0 and $t=t_1$. $(z,z')$ are natural complex coordinates for the
plumbing variety---see e.g.\ \cite{wolpert02}. Then as above there
is a complex $3g-3-n$ dimensional family of real-analytic
coordinates $s$ corresponding to quasiconformal maps which
preserve the complex coordinates near the nodes.
$(s,z,z',t_2,\cdots,t_n)$ form a real-analytic coordinate
neighborhood in $\overline{\mathcal M}_{g,1}$ so that for each
fixed $s$, $(z,z',t_2,\cdots,t_n)$ are complex-analytic
coordinates.

Let $\rho$ be a positive integer.  On any Riemann surface with
local coordinate $z$, a section $U$ of $K^\rho$ with a pole of
order $\rho$ at $p=\{z=0\}$ has a \emph{residue} $R^U_p$, which is
the the coefficient of the $\left(\frac{dz}z\right)^\rho$ term in
the Laurent series of $U$.  It is easy to check $R^U_p$ does not
depend on the complex coordinate $z$. Over a curve $C$ with $n$
nodes formed by pairing up $n$ pairs of points $(p_i,q_i)$ in
$\tilde C$, the space of \emph{regular $\rho$-differentials} over
$C$ is
$$\left\{U\in H^0(\tilde C, K_{\tilde C}^\rho\Pi_i([p_i]^\rho[q_i]^\rho)):
R^U_{p_i}=(-1)^\rho R^U_{q_i}\quad \forall \,i=1,\dots,n\right\}.
$$ In other
words these are sections of $K^\rho$ over $\tilde C$ with poles of
order $\rho$ allowed at $p_i,q_i$ so that the residues match up
appropriately on either side of each node. We are interested in
the case $\rho=3$ of regular 3-differentials and we denote the
sheaf of regular 3-differentials over $C$ as $K^{3,{\rm reg}}_C$.

Let $\mathcal{K}^3$ be the line bundle over $\overline{\mathcal
M}_{g,1}$ whose fiber over a pointed curve $(C,p)$ is the vector
space $K^{3,{\rm reg}}_C (p)$ the fiber of the regular
3-differentials over $C$ at $p$.  It straightforward that this
forms a holomorphic line bundle over $\overline{\mathcal
M}_{g,1}$, except possibly in the situation of a nodal curve $C$
in which the marked point $p$ is equal to a node.  Recall the
situation above: we separate the node and place a bubble
$\mathbb{CP}^1$ in between attached to the local normalization
$\tilde C ^{\rm loc}$ at points $r$ and $r'$. $p$ is any other
point in $\mathbb{CP}^1$.  Call this new curve $C^{\rm bubble}$.
Then any regular 3-differential on $C^{\rm bubble}$ has residues
$$R_{z=0}=-R_{r}\quad\mbox{and}\quad R_{z'=0}=-R_{r'}.$$
Now $\dim_\mathbb{C} H^0(\mathbb{CP}^1,K^3 [r]^3[r']^3)=1$, and
the residues satisfy $R_{r}=-R_{r'}$. Thus $R_{z=0}=-R_{z'=0}$. So
$\mathcal{K}^3$ at this point has a stalk naturally corresponding
to the residues at the node $p=\{z=z'=0\}$.

Now a local trivialization of $\mathcal{K}^3$ near the point in
$\overline{\mathcal{M}}_{g,1}$ corresponding to $(C,p)$
is given in terms of the coordinates $z,z'$ of the
plumbing variety.  For $z\neq0$, $z'\neq0$, let
$$T=\left(\frac{dz}z\right)^3= -\left(\frac{dz'}{z'}\right)^3,$$
and $T$ is well-defined except at the node. Then
$T$ extends over the node $\{z=z'=0\}$ to a local trivialization
so that at the node $T$ has residues $R^T_{z=0}=1$ and
$R^T_{z'=0}=-1$.  This is simply in terms of the coordinates
on the plumbing variety $z,z'$.  In general we take
\begin{equation} \label{trivialize}
T=c(s)\left(\frac{dz}z\right)^3= -\,c(s)\left(\frac{dz'}{z'}\right)^3,
\end{equation}
where the scalar $c(s)$ depends only real-analytically on $s$. We
may take $c(0)=1$.  This factor appears because the $s$ coordinates
on $\overline{\mathcal M}_{g,1}$ are only real-analytic with respect
to the underlying complex structure.

Now push forward the sheaf $\mathcal K^3$ by the map
$\mathcal F\!: \overline{\mathcal M}_{g,1} \to
\overline{\mathcal M}_g$
to form a sheaf $\mathcal S = K^3_{\overline{\mathcal M}_{g,1}/
\overline{\mathcal M}_g}$
over $\overline{\mathcal M}_g$.  In particular,
since the cohomology $H^i(C,K^{3,{\rm reg}})=\{0\}$ for all $i>0$,
and ${\rm dim}\,H^0(C,K^{3,{\rm reg}})=5g-5$, a theorem of
Grauert \cite{grauert60} shows that $\mathcal{S}$ is a (V-manifold)
vector bundle over $\overline{\mathcal M}_g$.  See Masur
\cite[Prop.\ 4.2]{masur76}
or Fay \cite{fay73} for details.

In particular, in a neighborhood of point corresponding to a Riemann
surface with nodes in $\overline{\mathcal M}_g$, there is a holomorphic
frame of regular 3-differentials.

\begin{prop} \label{U-neighborhood}
A basis for the analytic toplogy on the total space of $\mathcal{S}=
\mathcal{F}_*(\mathcal K^3)$
consists of neighborhood of $(C,U)$, where
$C$ is a curve with $n$ nodes and $U\in H^0(C,K^{3,{\rm reg}})$ consists
of pairs $(C_{s,t},U_{s,t})$ so that $(s,t)$ is close to zero in
$\co^{3g-3}$, and $U_{s,t}$ is close to $U$ in the following way:
Let the plumbing collars be represented by $\{z_iz'_i=t_i:|z_i|,|z'_i|<K\}$.
Outside the plumbing collars, we require
$$\left|\frac{U_{s,t}}{(dz^{\nu(s)})^3} - \frac{U}{dz^3}\right|<\epsilon$$
for $z^{\nu(s)}$ the local conformal coordinate on $C_{s,t}$
determined by the quasiconformal map for a local coordinate $z$ on
$C$.  Inside the plumbing collars $|z_i|,|z_i'|<K$, we require
$$\left|\frac{z_i^3}{dz_i^3}\left(U_{s,t}-U\right)\right|<\epsilon.$$
\end{prop}

\begin{proof}
This is proved in the above paragraphs.  In particular the last
statement follows from the trivialization
(\ref{trivialize}) of $\mathcal K^3$. We remark that since the cusp
coordinate $z_i^{\nu(s)}=e^{i\theta_s}z_i$ by Wolpert's Lemma 1.1 in
\cite{wolpert92b}, $ \frac{(dz_i^{\nu(s)})^3}{(z_i^{\nu(s)})^3} =
\frac{dz_i^3}{z_i^3}$ for $z_i$ inside the collar neighborhood.
\end{proof}

\begin{rem}
Fay \cite{fay73} and Yamada
\cite{yamada80} produce a more explicit asymptotic expansion
for a basis of regular 1-differentials on a degenerating family
of Riemann surfaces, and Masur
\cite{masur76} does the same for regular 2-differentials.  The same
techniques should apply to the present
case of regular 3-differentials as well.  Such a specific result
is not needed in this paper.
\end{rem}

We use this characterization to describe a degenerating family of
pairs $(C_\tau,U_\tau)$, which approach a pair $(C_0,U_0)$ of a
nodal curve $C_0$ and a regular 3-differential $U_0$ on $C_0$.  In
general, the parameter $\tau=\tau(s,t_1,\dots,t_n)$.  For
notational simplicity, we focus on the case $\tau=t_1=t$ the
plumbing parameter of the first node. Given a noded surface $C_0$
and a degenerating family $C_t$ determined by plumbing coordinates
$zz'=t$, a regular 3-differential $U_0$ on $C_0$ may be described
as a limit of holomorphic 3-differentials $U_t$ on $C_t$. Recall
(\ref{node-neigh}). Then
\begin{eqnarray*}
U_0|_{z'=0} &=& \sum_{m=-3}^\infty a_m z^m dz^3 \\
U_0|_{z=0} &=& \sum_{m=-3}^\infty b_m (z')^m (dz')^3
\end{eqnarray*}
where $b_{-3}=-a_{-3}$.  Then by Proposition \ref{U-neighborhood}
we may define on the collar neighborhood $\{z:|z|\in
[\frac{|t|}K,K]\} \subset C_t$
$$
U_t=\sum_{m=-\infty}^\infty a_m(t) z^m dz^3,
$$
where $a_m(t)$ is a continuous function of $t$ so that for
$m\ge-3$, $a_m(0)=a_m$, and for $m\le-3$,
$a_m(t)=-t^{-m-3}b_{-m-6}(t)$, where $b_n(t)$ is continuous in $t$
and $b_n(t)=b_n$ for $n\ge-3$.  (To make $U_t$ a holomorphically varying
family, we may choose $a_m$, $b_m$ to vary holomorphically in $t$
up to the matter of a branch cover of degree 2.  See e.g.\ Masur
\cite{masur76}.)

It is useful to compute in terms of more symmetric coordinates.
For some choice of branch of $\log t$, let $\ell$ be the
quasi-coordinate function given by
\begin{equation} \label{ell-coord}
\ell=\log z - \sfrac12 \log t = -\log z' +\sfrac12 \log t,
\end{equation}
and let $\mu=\rm{Re}\,\ell$.  Then we have
\begin{equation} \label{U-t-def}
 U_t= \left( -\sum_{n=1}^\infty b_{n-3}(t)t^{\frac{n}2}e^{-n\ell}
 + a_{-3}(t)
+ \sum_{n=1}^\infty a_{n-3}(t) t^{\frac{n}2}e^{n\ell} \right)
d\ell^3.
\end{equation}
Note that we may include the parameters $(s,t_2,\dots,t_n)$.  In
this case the coefficients $b_m$ and $a_m$ vary continuously in these
parameters as well.  Also, let $\ell_{s,t}=\log z^{\nu(s)} -\frac12\log t$
for $z^{\nu(s)}$ the complex coordinate on $C_{s,t}$ given by
the quasi-conformal map determined by $\nu(s)$.

\subsection{Holonomy of necks approaching a QH end} \label{qh-neck}
Recall that at a QH end, the residue $R=a_{-3}(0)\neq0$.  In this
case, we have on either side of the node, the model metric is
given by (\ref{flat-metric}), and which is the same for residue
$R$ and $-R$. These are flat cylindrical metrics of the same
radius which can simply be glued together in the plumbing
construction. For each $t$ small, we modify the cylindrical metric
to be
\begin{equation} \label{h-t}
h_t = 2^{\frac13}|a_{-3}(t)|^{\frac23}|d\ell|^2 .
\end{equation}
(Recall $|d\ell|^2=\frac{|dz|^2}{|z|^2}=\frac{|dz'|^2}{|z'|^2}$.)
Note that again, for simplicity, we repress the dependence of
$h_t$ on the other variables $(s,t_2,\dots,t_n)$, in which $h_t$
varies continuously. We may assume that the plumbing parameter $t$
satisfies $|t|< c^2$, where $c$ is the constant as in
(\ref{def-h}), so that the Ansatz metric is flat for $|z|,|z'|<c$.
(Note we can choose a uniform $c$ for a neighborhood of 0 in $t$
as long as $a_{-3}(t)$ is bounded away from 0.)

We will modify the barriers constructed in Subsection
\ref{find-sol} above to show that as $t\to0$, the solution $u$ to
(\ref{wang-eq}) will go to zero on the neck which is being pinched
to the node.  Although the cylindrical metrics fit together well,
the barriers must be modified.

Our barriers will be functions only of $\mu={\rm Re}\,\ell$. For
$S=S(\mu)$, use (\ref{U-t-def}) so that equation (\ref{operator})
becomes
\begin{equation}
\label{pert-L-eq}
L(S)=  2^{-\frac13}|a_{-3}(t)|^{-\frac23}S'' -2e^S + 2 e^{-2S}
\left[1+O(|t|^\frac12e^\mu)+O(|t|^\frac12e^{-\mu})\right].
\end{equation}
We need $L(S)\ge0$ for a lower barrier and $L(S)\le0$ for an upper
barrier.

Modify the barrier only in a neighborhood of the loop
$$\mathcal L_t =\{|z|=|z'|=|t|^\frac12\} = \{ \mu=0\},$$
where $\mu={\rm Re}\,\ell$.  We may do this because outside a
neighborhood $|z|,|z'|<K$, the original barriers suffice---outside
this neighborhood, $U_0$ changes to $U_t$ continuously, and all
the choices made in constructing the barriers in Subsection
\ref{find-sol} can easily be made to accommodate this small
perturbation.

Recall our upper and lower barriers in Subsection \ref{find-sol}
are of the form $\pm \beta|z|^{2\alpha}$, for $\beta\gg0$ and
$\alpha>0$ small. We may adjust these constants so that the same
$\beta$ and $\alpha$ are valid for both upper and lower barriers
on both sides of the puncture. On the plumbed surface choose our
upper barrier $S_t$ to be equal to
$$ S_t = \left\{ \begin{array}{l@{\quad\mbox{for}\quad}l}
\beta |t|^\alpha e^{2\alpha\mu} &
 \mu \in \left[1,\log K -\frac12 \log|t|\right] \\
I_t(\mu) &
 \mu \in [-1,1] \\
\beta |t|^\alpha e^{-2\alpha\mu} &
 \mu \in \left[-\log K +\frac12 \log|t|,-1\right]
\end{array}
\right.
$$
Notice that for the first and third lines of this
definition are respectively $\beta|z|^{2\alpha}$ and
$\beta|z'|^{2\alpha}$.  The middle part $I_t(\mu)$ is a $C^2$
interpolation between these two.  Explicitly, we may take
$I_t(\mu)$ to be the even fourth-order polynomial in $\mu$ so that
$S_t$ is $C^2$ at $\mu=\pm1$.  In other words, for
$Q=S_t(1)=S_t(-1)=\beta|t|^\alpha e^{2\alpha}$,
\begin{eqnarray*}
 I_t(\mu)& =& Q \left[ (1-\sfrac54\alpha+\sfrac12 \alpha^2)
  + (\sfrac32 \alpha - \alpha^2)\mu^2 +
   (-\sfrac14\alpha + \sfrac12 \alpha^2)\mu^4\right],  \\
I_t''(\mu) &=& \alpha Q\left[ (3-2\alpha)+(-3+6\alpha)\mu^2 \right].
\end{eqnarray*}
Therefore, we claim we can choose $\alpha$ independent of $t$ so
that for $t$ near 0, $I_t(\mu)>0$ and $L(I_t)\le0$ for
$\mu\in[-1,1]$. As in (\ref{L-eq}), rewrite (\ref{pert-L-eq}) as
\begin{eqnarray*}
L(I_t) &=& \left(2^{-\frac13}|a_{-3}(t)|^{-\frac23} I_t'' - 3I_t\right)
+ \left(2e^{-2I_t} - 2e^{I_t} +3I_t\right) \\
&& {}+ e^{-2I_t}\left[O(|t|^{\frac12}e^\mu)+O(|t|^\frac12
e^{-\mu})\right]
\end{eqnarray*}
$L(I_t)\le 0$  follows:  $2e^{-2I_t}-2e^{I_t}+3 I_t < 0$ for if
$I_t>0$. For $\alpha$ small, $Q$ dominates the perturbation terms
$O(|t|^\frac12 e^\mu)$ and $O(|t|^\frac12 e^{-\mu})$; therefore,
the first term dominates the last and $L(I_t)\le0$. This shows
that $S_t$ is an upper barrier for (\ref{wang-eq}) on the region
$\mu={\rm Re}\,\ell\in[-1,1]$.

For $|\mu|\in[1,\log K -\frac12 \log |t|]$, Proposition
\ref{U-neighborhood} shows $\|U_t\|_{h_t}^2\to \|U\|_h^2$
uniformly in in the plumbing variety coordinates $z,z'$ as
$t\to0$. By our choice of $h_t$ in (\ref{h-t}),
$$\|U_t\|_{h_t}^2 \le 2 + C|z| + C'|z'|$$
for constants $C,C'$.  Since this is true for all $t$, we may
choose $C,C'$ uniformly in $t$. Then the last term in (\ref{L-eq})
is dominated by the first for $\alpha$ small, $\beta$ large, and
$(z,z')$ close to 0, where these choices may be made independently
of $t$.  Then $L(S_t)\le 0$ for $|\mu|\in[1,\log
K-\frac12\log|t|]$.

Therefore, $S_t$ is an upper barrier for (\ref{wang-eq}).
Essentially the same arguments show that near $\{|\mu|\le\log K
-\frac12\log|t|\}$, $-S_t$ forms a lower barrier for solutions to
(\ref{wang-eq}) with data $(C_t^{\rm reg},U_t,h_t)$, and as in
Subsection \ref{find-sol}, there is a constant $B_t$ so that a
lower barrier is equal to $-S_t$ inside the plumbing collars and
$B_t$ outside.  $B_t$ depends continuously on $t$.

Again following Subsection \ref{find-sol}, the maximum principle
says that there is a bounded solution $u_t$ to (\ref{wang-eq}) on
$C_t^{\rm reg}$ so that $|u_t| \le S_t$ for $|z|,|z'|<K$.  Note
that as $t\to0$, $S_t\to0$ on a neighborhood of the loop $\mathcal
L_t$. The geometry is still uniformly bounded in this case; so we
still have the bound on $u_t$ and $\nabla u_t$ as in Proposition
\ref{bound-prop}. Then let $\mathbf{A}_t$ be the holonomy around
the loop $\mathcal L_t$ with respect to the frame $\{f,f_w,f_{\bar
w}\}$ as in (\ref{x-hol}). Orient $\mathcal L_t$ counterclockwise
in the $z$ coordinate. Then Lemma \ref{A-lemma} still holds and we
have as in the proof of Proposition \ref{hyperbolic-prop}

\begin{prop}
As $t\to0$, the eigenvalues of the holonomy along $\mathcal L_t$
approach $e^{2\pi\lambda_i}$ for $\lambda_i$ the roots, with
multiplicity,  of formula (\ref{char-poly}).
\end{prop}

In terms of the other parameters $(s,t_2,\dots,t_n)$, the Ansatz
metric $h_{s,t}$ and the upper and lower barriers vary
continuously. One thing to note is that for $s$ small, Wolpert's
Beltrami differential $\nu(s)$ is close to 0 and supported away
from the collar neighborhoods.  Therefore, the complex structure
and hyperbolic metric on $C^{\rm reg}_s$ are close to that on
$C^{\rm reg}_0$. Then on $C^{\rm reg}_{s,t}$, we take $h_{s,t}$ to
be the hyperbolic metric on $C^{\rm reg}_{s,0}$ outside the collar
neighborhoods and modify the metrics $h_t$ as in
(\ref{def-h},\ref{h-t}) inside the collar neighborhoods only.

\subsection{Vertical twist parameters of a family} \label{vert-twist-fam}
Recall that if ${\rm Re}\,R>0$, then the vertical twist parameter
of the end is $+\infty$.  Theorem \ref{families} part
(\ref{twist-fam}) will follow from Proposition
\ref{holo-matrix-param} above.  First note that in a neighborhood
in the total space of $\mathcal K^3$, by the previous subsection,
there are uniform bounds on $\|U_\tau\|$ on each Riemann surface
$C_\tau^{\rm reg}$ for $\tau$ near 0.  Also, there are uniform
barriers of the form $\beta |z_i|^{2\alpha}$ and $\beta
|z'_i|^{2\alpha}$ in each plumbing collar, and uniform constant
barriers in on the rest of $C_\tau$. Then we may bootstrap as in Lemma
\ref{gradient-bound} to find uniform estimates on the third
covariant derivatives of $u_\tau$ the solution to (\ref{wang-eq})
for data $(C_\tau^{\rm reg},U_\tau,h_\tau)$. Therefore, by
Ascoli-Arzela, a subsequence of $u_\tau$ converges in $C^2$
(determined by covariant derivatives of the metric) to a solution
$\tilde u$ of (\ref{wang-eq}) for data $(C^{\rm reg},U,h)$.  By
Proposition \ref{unique}, $\tilde u=u$ and thus there are $C^2$
estimates for $u_\tau$ approaching $u$.  In particular, by
(\ref{y-eq}), the perturbation matrix $\mathbf{B}^{\rm pert}_\tau$
from (\ref{y-dir}) varies continuously in $\tau$ in the plumbing
collar.

More specifically, in terms of the $\ell$ coordinate in a collar
neighborhood as in the previous subsection, $\mathbf{B}^{\rm
pert}_\tau$ is uniformly continuous in $\tau$ and satisfies
$$|\mathbf{B}^{\rm pert}_\tau|\le \beta |t|^\alpha e^{2\alpha
\mu} \quad \mbox{for } \mu\in[1,\log K-\sfrac12\log|t|],
$$ where $t=t(\tau)$. In terms of the $y$ coordinate,
$$|\mathbf{B}^{\rm pert}_\tau|\le \beta e^{-2\alpha y}
\quad \mbox{for } y\in[-\log K,-1-\sfrac12\log|t|].
$$ We have similar bounds in the $\tilde y$ coordinates
(\ref{ytilde-def}). Now for $t$ small, replace $\mathbf{B}^{\rm
pert}_\tau$ by
$$\mathbf{B}^{\rm cutoff}_\tau= \left\{
\begin{array}{r@{\mbox{ for }}l}
\mathbf{B}^{\rm pert}_\tau & y\le -1-\frac12\log|t| \\
f(x,y,\tau) & y>-1-\frac12\log|t|
\end{array} \right.
$$
so that $\mathbf{B}^{\rm cutoff}_\tau$ satisfies the hypotheses of
Lemma \ref{pert-sol}. Then for $t=t(\tau)$ small, the solution to
(\ref{y-dir}) satisfies the asymptotic bounds (\ref{z-sol}), and
moreover the $o(e^{\mu_i\tilde y})$ terms are uniform in $\tau$ as
$t(\tau)\to0$.  Therefore, as in Subsection
\ref{full-triangle-sec} above, we may fix a large initial value of
$y$, say $y_0$, independent of $\tau$, and the solutions to the
equations (using $\mathbf{B}^{\rm cutoff}$) in the two directions
determined by $\iota$ and $\hat\iota$ give solutions which
approach the geodesic line segments $G^{+0}$ and $G^{0-}$
respectively. (Notice that in these cases
by (\ref{a-hat}) $\tilde{\mathbf B}$ and
$\mathbf{B}^{\rm cutoff}$
satisfy the hypotheses of Proposition \ref{app-prop}.)
Now for the actual $\rp^2$ structure determined by $(C_\tau^{\rm
reg},U_\tau)$, we must take our original perturbation term
$\mathbf{B}^{\rm pert}$. But for $y\le -1-\frac12\log|t|$, the
solutions are the same by uniqueness of solutions to ODEs.  Since
$y_0$ is fixed and the bound $-1-\frac12\log|t| \to\infty$ as
$t\to0$, there are, as $t\to0$, two regions of the developing map
which approach the geodesic line segments $G^{+0}$ and $G^{0-}$
respectively.
The only way this can happen is if the vertical twist parameter is
approaching $+\infty$ as $\tau\to0$.

In the case of a node with ${\rm Re}R_i<0$, there is not as much
information concerning the holonomy maps, but we do know that the
residue of the other half of the node has positive real part, and
from the point of view of the residue with ${\rm Re}R_i>0$, the
vertical twist parameter approaches $+\infty$.  Therefore, the
vertical twist parameter from the opposite point of view shrinks
to $-\infty$.

\subsection{Holonomy of necks approaching a parabolic end}
Recall from Section \ref{finding-sol} that the Ansatz and the
barrier for a parabolic end (one for which the residue $R_i=0$) is
of quite a different form.  The Ansatz metric is the hyperbolic
metric, and the upper and lower barriers near the puncture are
both constants.  We are free to make these constants larger in
norm by the methods above, and thus we can assure that the
barriers on either side of the node are the same constant, and
thus patch together naturally as the node is plumbed.

The metric, on the other hand, must be smoothed across the plumbed
node.  It is crucial that the curvature still be negative and
bounded away from 0 and $-\infty.$ For this purpose, we recall the
plumbing metric in Wolpert \cite{wolpert90}.  See also
Wolf-Wolpert \cite{wolf-wolpert92}.

\begin{prop} \cite{wolpert90}  \label{kappa-bound}
Let $C_0$ be a stable nodal Riemann surface, and let $C_t$ be the
plumbed surface as above. On the plumbing collar $\{zz'=t\}$ for
$|z|,|z'|\le K$, there is a metric $h_t$ which is equal to the
hyperbolic metric on $C_0$ outside the plumbing collar(s), and
satisfies
\begin{equation} \label{graft-metric}
h_t = \left\{
\begin{array}{l@{\mbox{  for  }}l@{\mbox{  i.e.,  }}l}
\frac{4|dz|^2}{|z|^2(\log|z|^2)^2} & |z|\in [\frac{2K}3,K] &
|z'|\in[ \frac{|t|}K,\frac{3|t|}{2K} ] \\

\left(\frac\pi{\log|t|} \csc(\frac{\pi\log|z|}{\log|t|})
\frac{|dz|}{|z|} \right)^2 &
 |z| \in [\frac{3|t|}K ,\frac{K}3] & |z'| \in [\frac{3|t|}K ,\frac{K}3] \\

\frac{4|dz'|^2}{|z'|^2(\log|z'|^2)^2} & |z|\in
[\frac{|t|}K,\frac{3|t|}{2K} ] &
 |z'|\in[\frac{2K}3,K] \\
\end{array} \right.
\end{equation}
and smoothly interpolated for $|z|,|z'|\in[\frac{K}3,\frac{2K}3]$.
The curvature $\kappa_t$ on $(C_t,h_t)$ satisfies
$$ \|\kappa_t+1\|_{C^0} \le \gamma|t|^\delta$$
for uniform constants $\gamma,\delta>0$.
\end{prop}

\begin{rem}
For our purposes, we may take the model grafting in Section 3.4.MG of
Wolpert \cite{wolpert90}.  We do not need to use the more
complicated (and accurate) grafting procedures in Section 3.4.CG of
\cite{wolpert90} or in \cite{wolf-wolpert92}.
\end{rem}

\begin{lem} \label{u-bound-lemma}
There is a constant $\delta'$ so that if $|t|\le\delta'$ and
$\|U_t\|_{h_t} \le \mathcal K$ for a constant $\mathcal K$
independent of $t$, then there is a constant $\mathcal K'$
independent of $t$ so that solution $u_t$ to (\ref{wang-eq})
satisfies $|u_t|\le \mathcal K'$.
\end{lem}

\begin{proof}
In Subsection \ref{find-sol}, the upper and lower barriers for
(\ref{wang-eq}) on the curve $C_0$ is equal to a constant in the
neighborhood of the node in question (note that the constants may
be adjusted on either side of the node to be equal).  A constant
$M_t>0$ is an upper barrier of (\ref{wang-eq}) if
\begin{eqnarray*}
4e^{-2M_t}\|U_t\|_{h_t}^2 -2 e^{M_t}  - 2 \kappa_t &\le& 0 \quad
\mbox{i.e.,} \\
4 \|U_t\|_{h_t}^2 - 2 E_t^3 - 2\kappa_t E_t^2 &\le& 0,
\end{eqnarray*}
for $E_t=e^{M_t}$.  It is easy to see that $E_t$ may be chosen
independently of $t$ given that $\|U_t\|_{h_t}$ is bounded
independently of $t$, and (by Proposition \ref{kappa-bound})
$\kappa_t$ satisfies $-k'\le\kappa_t\le-k$ for positive constants
$k,k'$.  Similar considerations apply for the lower barrier.
\end{proof}

\begin{prop}
Let $\mathbf{A}_t$ be the holonomy around the loop $\mathcal L_t=
\{|z|=|z'|=|t|^\frac12\}$ with respect to the frame
$\{f,f_w,f_{\bar w}\}$ as in (\ref{x-hol}). Orient $\mathcal L_t$
counterclockwise in the $z$ coordinate.  If there are uniform
positive constants $\delta', C$ so that for $|t|\le\delta'$,
\begin{equation} \label{U-t-assumption}
\sup_{|z|\in[\frac{K}{|t|},K]}
\left|\frac{z^3(U_t-U_0)}{dz^3}\right| \left|\log|t|\right|^3  \le C,
\end{equation}
then $\mathbf{A}_t$ is continuous in $t$ and
$$\lim_{t\to0}\mathbf{A}_t = \left(\begin{array}{ccc}
0&1&1 \\ 0&0&0 \\ 0&0&0 \end{array}\right).$$
\end{prop}

\begin{proof}
Recall
$$ \mathbf{A}_t = \left( \begin{array}{ccc}
0&1&1\\
\frac12e^{\psi} & \psi_w& Ue^{-\psi} \\
\frac12e^{\psi} & \bar{U}e^{-\psi} & \psi_{\bar w}
\end{array} \right).$$
for $z=e^{iw}$, $\psi=\phi+u$, $h_t=e^\phi|dw|^2$, and $w=x+iy$.
Note we suppress the dependence on $t$. The continuity in $t$
follows from the fact that $u_t$ and its derivatives vary
continuously in $t$ by standard elliptic regularity arguments, as
in the first paragraph of Subsection \ref{vert-twist-fam} above.
(In particular, it is easy to check that the metric
(\ref{graft-metric}) has bounded geometry on a uniformly large
neighborhood of $\mathcal L_t$.) Also near $\mathcal L_t$,
$$ \phi = 2\log\left(\frac{-\pi}{\log|t|}\right) - 2 \log \sin
\left(\frac{-\pi y}{\log|t|} \right).
$$ On $\mathcal L_t$, $y=-\frac12\log|t|$, and so
$$ e^\phi=\left(\frac\pi{\log|t|}\right)^2, \quad \mbox{and}
\quad \phi_w=0.
$$
Moreover, by Lemma \ref{u-bound-lemma}, $|u|\le\mathcal K'$.  We
still have Lemma \ref{gradient-bound}, which shows
$$|u_w|\le
\mathcal K '' e^{\frac\phi2} = \mathcal K''
\frac\pi{\left|\log|t|\right|}
$$
Finally, the assumption (\ref{U-t-assumption}) on $U$, the fact
that $U_0$ has residue 0, and the uniform bound on $u$ shows as
$t\to0$, $Ue^{-\psi}\to0$ on $\mathcal L_t$. Similarly, all the
other entries in the second and third rows of $\mathcal{A}_t$ go
to zero.
\end{proof}

Then as above in Subsection \ref{main-holo-calc}, the eigenvalues of
the holonomy around the loop $\mathcal L_t$ all approach 1 as $t\to0$.

In terms of more general paths $(C_\tau,U_\tau)$ in $\mathcal S$,
calculating the limiting holonomy depends on having uniform $C^0$
estimates on $u$ independent of $\tau$.
The model metric on $C_\tau$ is simply the hyperbolic metric on
$C_{s(\tau),0}$ outside the collar neighborhoods and may be modified
as above in each collar neighborhood.  Note that the complex structure
and hyperbolic metrics on $C_{s(\tau),0}$ vary continuously as
in Wolpert's Lemma \cite{wolpert92b}.
In a neighborhood
of a singular curve $C_0$ with $n$ nodes, consider V-manifold coordinates
$(t_1,\dots,t_n,s)$ in $\overline{\mathcal{M}}_g$ near $C_0$.  Let
$U_0$ have residue 0 at $k$ of the $n$ nodes.  Without loss of generality,
assume these nodes correspond to the plumbing parameters
$(t_1,\dots,t_k)$.  Then

\begin{prop}
For $(C_\tau,U_\tau)$ a continuous path in $\mathcal S$, if in addition
there are uniform positive constants $\delta,C$ so that
for $|\tau|\le\delta$, $U_\tau$ satisfies
\begin{equation} \label{U-tau-bound}
\sup_{|z_i|,|z_i'|\in[\frac{K}{|t_i|},K]}
\left|\frac{z_i^3(U_\tau-U_0)}{dz_i^3}\right|
\left|\log|t_i|\right|^3  \le C
\end{equation}
for $t_i=t_i(\tau)$ and for all $i\in\{1,\dots,k\}$, then the eigenvalues
of the holonomy around each neck are continuous in $\tau$.
\end{prop}

Given a holomorphic frame $\{\Psi^1,\dots,\Psi^{5g-5}\}$ of the
vector bundle $\mathcal S\to\overline{\mathcal M}_g$, we have the
following

\begin{cor}
Write $U_\tau = a_j(\tau)\Psi^j$, where $a=(a_j)\in\co^{5g-5}$ and
$\Psi^j$ represents the element of $H^0(C_\tau,K^{3,{\rm reg}})$
corresponding to the frame. Then if there is a uniform $C$ so that
$$ |[a(\tau)-a(0)] (\log|t_i|)^3|\le C \quad \mbox{for }i=1,\dots,k,$$
then the eigenvalues of the holonomy around each neck are continuous
in $\tau$.
\end{cor}

\subsection{Results} We record the results of the previous
subsections in

\begin{thm} \label{families}
Consider a continuous path of pairs $(C_\tau,U_\tau)$, where the
possibly nodal curve $C_\tau$ represents a point in
$\overline{\mathcal{M}}_g$ and $U_\tau$ is a holomorphic section
of $K^{3,{\rm reg}}_{C_\tau}$ so that $C_0$ is a nodal curve with
$n$ nodes.  For each node, pick one side from which to measure the
residue.  For any curve $C_\tau$ which approximates $C_0$ by
pinching a neck to form the node, this amounts to choosing an
orientation for any loop around that neck. Then $U_0$ is a cubic
differential with residue $R_i$ for each node.
\begin{enumerate}
\item If all the residues $R_i\neq0$, then the eigenvalues of the
holonomy around each neck which is pinched to the node
continuously approach the eigenvalues of the holonomy around the
punctures of the complete Riemann surface $C_0^{\rm reg}$ as in
Table \ref{holonomy-table}.  The same is true if we have some
residues $R_1,\dots,R_k=0$ as long as $U_\tau$ satisfies the
addition set of bounds (\ref{U-tau-bound}).

\item \label{twist-fam} Still assume (\ref{U-tau-bound}) for all
nodes with 0 residue. Consider a node whose residue $R_i$
satisfies ${\rm Re}\,R_i\neq0$. Then the vertical twist parameter
along this neck $N$ approaches $\pm\infty$, the sign agreeing with
the sign of the ${\rm Re}\,R_i$.  In fact, if ${\rm Re}\,R_i>0$,
there is a continuous path of points $p_\tau\in \overline{\mathcal
M}_{g,1}$ so that $p_\tau\in C_\tau$ and $p_\tau$ avoids all
nodes, and a continuously varying choice of $\rp^2$ coordinate
chart near $p_\tau$ in the $\rp^2$ surface $S_\tau$ determined by
$(C_\tau^{\rm reg},U_\tau)$.  The holonomy matrix with respect to
these $\rp^2$ coordinates of the neck $N$ has fixed points
Fix$^0_\tau$, Fix$^+_\tau$ and Fix$^-_\tau$ which vary
continuously with $\tau$.  If we fix coordinates on $\rp^2$ so
that Fix$^0_\tau$, Fix$^+_\tau$ and Fix$^-_\tau$ are fixed, then
the image $\Omega_\tau$ of the developing map satisfies
$$\lim_{\tau\to0}\Omega\tau \supset T,$$
where $T$ is the principal triangle whose vertices are the fixed
points.  (In other words for all points $q\in T$, there is a constant
$\delta$ so that if $|\tau|\le\delta$, then $q\in\Omega_\tau$.)
\end{enumerate}
\end{thm}

\begin{rem}
We expect that the technical restrictions on the continuous paths
$(C_\tau,U_\tau)$ in the case $R_i=0$ can be removed.
\end{rem}

\subsection{The $\rp^2$ structure on a degenerating neck}
By the definition of regular 3-differentials, it is worthwhile to
compare the $\rp^2$ holonomy of two punctures in $\Sigma=C^{\rm
reg}$ equipped with cubic differentials with residues $R$ and
$-R$, as these will naturally be identified in the nodal curve
$C$. Recall that the eigenvalues of the holonomy are given by
$e^{2\pi\lambda_i}$, where $\lambda_i$ are the roots of
(\ref{char-poly}).  If we replace $R$ by $-R$ in
(\ref{char-poly}), then the roots $\lambda_i$ become $-\lambda_i$.
In terms of the holonomy matrix, at least in the hyperbolic case,
the holonomy matrix satisfies $H_{-R}=H_R^{-1}$. We may think of
this as the same holonomy viewed from opposite orientations, which
is natural: the holonomy is given in terms of loops which go
counterclockwise around each puncture, and if want to glue two
such punctures together, the
two loops will be oriented in opposite directions.

In the case of parabolic and quasi-hyperbolic holonomies, the
holonomy is the only invariant we have of the end, but there is also
the vertical twist parameter for hyperbolic ends.
Assume Re$\,R>0$.  Then the vertical twist parameter for this end
is $+\infty$, while the vertical twist parameter with the
corresponding end with residue $-R$ will be $-\infty$.  Again this
is to be expected.  We may imagine a family of surfaces
degenerating in a way that their vertical twist parameters become
infinite.  Then measured from one side, (vertical twist parameter
going to $+\infty$) the piece of the developing map glued along
the principal geodesic is becoming larger and larger until it
becomes the entire principal triangle.  From the other side (for
which the vertical twist parameter goes to $-\infty$), the glued
piece of the developing map becomes smaller and smaller until it
vanishes and the boundary is simply the principal geodesic. See
Section \ref{goldman-coor}.

\appendix
\section{Linear almost constant-coefficient systems with parameters}
\label{app}

\begin{prop} \label{app-prop}
Consider a system of linear differential equations
\begin{equation} \label{app-eq}
 \partial_y X(s,y) =(c(s)B+R(s,y))X(s,y),
\end{equation}
where $y\ge T$, $s$ is a set of parameters,
$B=\mathbf{D}(\mu_1,\dots,\mu_n)$ is a constant $n\times n$
diagonal matrix, the scalar factor $c(s)$ is continuous, and the
matrix entries of the error term $R(s,y)$ are smooth and $L^1$ in
$y$ and vary in $s$ so that for all $i,j$, the  map $s \mapsto
R^j_i(s,y)$ is continuous to $L^1([T,\infty))$. 

Then there exist $n$ linearly independent solutions $X^{(k)}(s,y)$
to (\ref{app-eq}) so that for $\{v_k\}$ the standard basis on
$\re^n$,
$$X^{(k)}(s,y)=e^{c(s)\mu_k y}v_k + e^{c(s)\mu_k y} b(s,y),$$
and $X^{(k)}(s,y)$ is continuous in $(s,y)$,
$\lim_{y\to\infty}|b(s,y)|=0$ uniformly in $s$ for $s$ in a
bounded region.
\end{prop}

\begin{rem}
If in addition (as we have above), there are uniform positive
constants $C,\gamma$ so that $|R^j_i(s,y)|\le Ce^{-\gamma y}$,
then we can replace $c(s)B$ by a continuous family $B(s)$ of
diagonal matrices, and moreover a more precise error bound on
$X^{(k)}-e^{B(s)\mu_k y}v_k$ holds. Since we do not need this
better result, we do not prove it here.
\end{rem}

\begin{proof}
We follow the treatment in Levinson \cite[pp.\
115--117]{levinson48}. The only additional thing to prove is the
continuity of solutions in $s$. Assume $\mu_i$ the eigenvalues of
$B$ are arranged so that
$${\rm Re}\,\mu_1 \ge {\rm Re}\,\mu_2 \ge \cdots \ge {\rm Re}\,\mu_n.$$
Fix $k\in\{1,\dots,n\}$. Choose $q=q(k)$ so that
$${\rm Re}\,\mu_k={\rm Re}\,\mu_q>{\rm Re}\,\mu_{q+1}$$
(or $q=n$ if this is impossible). Define $X_i^0 =
\delta_{ik}e^{c(s)\mu_k y}$, and define $X_i^m$ recursively by
\begin{equation*}
\begin{array}{rclr}
X^{m+1}_i(s,y) &=& \D \delta_{ik}e^{c(s)\mu_k y}
                   - \int_y^\infty e^{\mu_i(y-\sigma)}
           R_i^j(s,\sigma) X^m_j(s,\sigma) \,d\sigma & (i\le q), \\
X^{m+1}_i(s,y) &=& \D \int_a^y e^{\mu_i(y-\sigma)}
           R_i^j(s,\sigma) X^m_j(s,\sigma) \,d\sigma & (i>q ).
\end{array}
\end{equation*}
The index $j$ is summed from $1$ to $n$, but no sum is taken over
$i$. If for some $m$, $X_i^m=X_i^{m+1}$, then this $X_i^m$ solves
(\ref{app-eq}) as long as the integrals involved converge
absolutely.

Now it is clear that $X_i^0$ is continuous in $(s,y)$.  By
induction, the same is true for $X_i^m$ for all $m$.  As in
\cite{levinson48},
$$|X^{m+1}_i - X^m_i|\le 2^{-m} \left|e^{c(s)\mu_k y}\right|,$$
and so the series (for $X^{-1}_i=0$)
$$ X_i=\sum_{m=-1}^\infty (X^{m+1}_i - X^m_i)$$
is majorized by a geometric series. Also $X=(X_i)$ solves
(\ref{app-eq}). In particular, $X^m$ converges locally uniformly
to $X$ in $(s,y)$ and $X$ is continuous in $(s,y)$.  This $X$ is
the $X^{(k)}$ in the proposition.

The bound on the error term in \cite{levinson48} is of the form
$$ |b_i(s,y)|\le e^{-\epsilon y} |R|_{L^1} + 2\sum_j
\int_{\frac{y}2}^\infty |R^j_i(s,\sigma)|\,d\sigma
$$ for a uniform positive constant $\epsilon$.
This shows the required uniform bound on $|b(s,y)|$. The rest of
the proposition follows as in \cite{levinson48}.
\end{proof}

\bibliographystyle{abbrv}
\bibliography{thesis}

\begin{thebibliography}{10}

\bibitem{baily57}
W.~L. Baily.
\newblock On the imbedding of {V}-manifolds in projective space.
\newblock {\em American Journal of Mathematics}, 79:403--430, 1957.

\bibitem{calabi72}
E.~Calabi.
\newblock Complete affine hyperspheres {I}.
\newblock {\em Instituto Nazionale di Alta Matematica Symposia Mathematica},
  10:19--38, 1972.

\bibitem{cheng-yau75}
S.-Y. Cheng and S.-T. Yau.
\newblock Differential equations on {R}iemannian manifolds and their geometric
  applications.
\newblock {\em Communications on Pure and Applied Mathematics}, 28:333--354,
  1975.

\bibitem{cheng-yau77}
S.-Y. Cheng and S.-T. Yau.
\newblock On the regularity of the {M}onge-{A}mp{\`e}re equation
  $\det((\partial^2u/\partial x^i\partial x^j))={F}(x,u)$.
\newblock {\em Communications on Pure and Applied Mathematics}, 30:41--68,
  1977.

\bibitem{cheng-yau86}
S.-Y. Cheng and S.-T. Yau.
\newblock Complete affine hyperspheres. part {I}. {T}he completeness of affine
  metrics.
\newblock {\em Communications on Pure and Applied Mathematics}, 39(6):839--866,
  1986.

\bibitem{choi94a}
S.~Choi.
\newblock Convex decompositions of real projective surfaces. {I}: $\pi$-annuli
  and convexity.
\newblock {\em Journal of Differential Geometry}, 40(1):165--208, 1994.

\bibitem{choi94b}
S.~Choi.
\newblock Convex decompositions of real projective surfaces. {II}: Admissible
  decompositions.
\newblock {\em Journal of Differential Geometry}, 40(2):239--283, 1994.

\bibitem{choi96}
S.~Choi.
\newblock Convex decompositions of real projective surfaces. {III}: For closed
  and nonorientable surfaces.
\newblock {\em Journal of the Korean Mathematical Society}, 33(4):1139--1171,
  1996.

\bibitem{choi-goldman93}
S.~Choi and W.~M. Goldman.
\newblock Convex real projective structures on closed surfaces are closed.
\newblock {\em Proceedings of the American Mathematical Society},
  118(2):657--661, 1993.

\bibitem{choi-goldman97}
S.~Choi and W.~M. Goldman.
\newblock The classification of real projective structures on compact surfaces.
\newblock {\em Bulletin (New Series) of the American Mathematical Society},
  34(2):161--171, 1997.

\bibitem{d-goldman}
M.-R. Darvishzadeh and W.~M. Goldman.
\newblock Deformation spaces of convex real projective and hyperbolic affine
  structures.
\newblock {\em Journal of the Korean Mathematical Society}, 33:625--638, 1996.

\bibitem{dunkel}
O.~Dunkel.
\newblock Regular singular ponts of a system of homogeneous linear differential
  equations of the first order.
\newblock {\em Proceedings of the American Academy of Arts and Sciences},
  38:341--370, 1902--03.

\bibitem{fay73}
J.~D. Fay.
\newblock {\em Theta Functions on {R}iemann Surfaces}, volume 352 of {\em
  Lecture Notes in Mathematics}.
\newblock Springer-Verlag, 1973.

\bibitem{gilbarg-trudinger}
D.~Gilbarg and N.~Trudinger.
\newblock {\em Elliptic Partial Differential Equations of Second Order}.
\newblock Springer-Verlag, 1983.

\bibitem{goldman84}
W.~M. Goldman.
\newblock The symplectic nature of the fundamental groups of surfaces.
\newblock {\em Advances in Mathematics}, 54:200--225, 1984.

\bibitem{goldman90a}
W.~M. Goldman.
\newblock Convex real projective structures on compact surfaces.
\newblock {\em Journal of Differential Geometry}, 31:791--845, 1990.

\bibitem{goldman90b}
W.~M. Goldman.
\newblock The symplectic geometry of affine connections on surfaces.
\newblock {\em Journal f{\"u}r die reine und agnewandte Mathematik},
  407:126--159, 1990.

\bibitem{grauert60}
H.~Grauert.
\newblock {E}in {T}heorem der analytischen {G}arbentheorie und die
  {M}odulr{\"a}ume komplexer {S}trukturen.
\newblock {\em Pub. Math. IHES}, 5:233--292, 1960.

\bibitem{hartman}
P.~Hartman.
\newblock {\em Ordinary Differential Equations}.
\newblock Wiley, 1964.

\bibitem{hartman-wintner55}
P.~Hartman and A.~Wintner.
\newblock Asymptotic integrations of linear differential equations.
\newblock {\em American Journal of Mathematics}, 77:45--86, 1955.

\bibitem{hitchin87}
N.~J. Hitchin.
\newblock The self-duality equations on a {R}iemann surface.
\newblock {\em Proceedings of the London Mathematical Society}, 55(3):59--126,
  1987.

\bibitem{hitchin92}
N.~J. Hitchin.
\newblock Lie groups and {T}eichm\"uller space.
\newblock {\em Topology}, 31(3):449--473, 1992.

\bibitem{kim99}
H.~C. Kim.
\newblock The symplectic global coordinates on the moduli space of real
  projective structures.
\newblock {\em Journal of Differential Geometry}, 53(2):359--401, 1999.

\bibitem{labourie97}
F.~Labourie.
\newblock in Proceedings of the GARC Conference in Differential Geometry, Seoul
  National University, Fall 1997, 1997.

\bibitem{levinson48}
N.~Levinson.
\newblock The asymptotic nature of solutions of linear systems of differential
  equations.
\newblock {\em Duke Mathematical Journal}, 15:111--126, 1948.

\bibitem{loftin01}
J.~C. Loftin.
\newblock Affine spheres and convex $\rp^n$ manifolds.
\newblock {\em American Journal of Mathematics}, 123(2):255--274, 2001.

\bibitem{masur76}
H.~Masur.
\newblock The extension of the {W}eil-{P}etersson metric to the boundary of
  {T}eichm{\"{u}}ller space.
\newblock {\em Duke Mathematical Journal}, 43(3):623--635, 1976.

\bibitem{nomizu-sasaki}
K.~Nomizu and T.~Sasaki.
\newblock {\em Affine Differential Geometry: Geometry of Affine Immersions}.
\newblock Cambridge University Press, 1994.

\bibitem{wang91}
C.-P. Wang.
\newblock Some examples of complete hyperbolic affine $2$-spheres in
  $\mathbb{R}^3$.
\newblock In {\em Global Differential Geometry and Global Analysis}, volume
  1481 of {\em Lecture Notes in Mathematics}, pages 272--280. Springer-Verlag,
  1991.

\bibitem{wolf91}
M.~Wolf.
\newblock Infinite energy harmonic maps and degeneration of hyperbolic surfaces
  in moduli space.
\newblock {\em Journal of Differential Geometry}, 33:487--539, 1991.

\bibitem{wolf-wolpert92}
M.~Wolf and S.~Wolpert.
\newblock Real analytic structures on the moduli space of curves.
\newblock {\em American Journal of Mathematics}, 114:1079--1102, 1992.

\bibitem{wolpert02}
S.~Wolpert.
\newblock Geometry of the {W}eil-{P}etersson completion of {T}eichm{\"{u}}ller
  space.
\newblock preprint, June 2002.

\bibitem{wolpert85}
S.~Wolpert.
\newblock On the {W}eil-{P}etersson geometry of the moduli space of curves.
\newblock {\em American Journal of Mathematics}, 107:969--997, 1985.

\bibitem{wolpert90}
S.~Wolpert.
\newblock The hyperbolic metric and the geometry of the universal curve.
\newblock {\em Journal of Differential Geometry}, 31(2):417--472, 1990.

\bibitem{wolpert92b}
S.~Wolpert.
\newblock Spectral limits for hyperbolic surfaces {I}{I}.
\newblock {\em Inventiones Mathematicae}, 108(1):91--129, 1992.

\bibitem{yamada80}
A.~Yamada.
\newblock Precise variational formulas for {A}belian differentials.
\newblock {\em Kodai Mathematical Journal}, 3:114--143, 1980.

\end{thebibliography}

\end{document}